\documentclass{article}
\usepackage{algorithmicx,algorithm}
\usepackage[pagewise]{lineno}\linenumbers
\usepackage{algpseudocode}
\usepackage{amsmath}
\usepackage{amssymb}
\usepackage{amsthm}
\usepackage{appendix}
\usepackage{array}
\usepackage{arydshln}
\usepackage{bbm}
\usepackage{bm}
\usepackage{color}
\usepackage{float}
\usepackage{geometry}
\usepackage{graphicx}
\usepackage{indentfirst}
\usepackage[utf8]{inputenc}
\usepackage{makecell,rotating,multirow,diagbox}
\usepackage{mathtools}
\usepackage{pgf,tikz,pgfplots}
\usetikzlibrary{positioning,arrows.meta}
\usepackage{slashbox}
\usepackage{subfig}
\usepackage{epstopdf}
\usepackage{xcolor}
\usepackage{setspace}
\usepackage{stackengine}
\usepackage{ulem}
\usepackage{enumitem}

\makeatletter 
\@addtoreset{equation}{section} 
\makeatother 
\renewcommand{\theequation}{\arabic{section}.\arabic{equation}}

\def \pt {\partial}
\def \eq$#1${\begin{equation}#1\end{equation}}
\newcommand{\dd}{\textup{d}}
\newcommand{\blue}[1]{{\color{black}#1}}

\newtheorem{theorem}{Theorem}[section]

\newtheorem{remark}[theorem]{Remark}

\newtheorem{proposition}[theorem]{Proposition}



\def \bpsi {\bm{\psi}}
\def \balpha {\bm{\alpha}}
\def \bbeta {\bm{\beta}}
\def \bu {\bm{u}}

\def \bomega {\bar{\omega}}
\def \bV {\bar{V}}
\def \bb {\bm{b}}
\def \be {\bm{e}}
\def \bmu {\bm{\mu}}
\def \bxi {\bm{\xi}}
\def \etab {\bm{\eta}}
\def \bc {\bm{c}}
\def \bs {\bm{s}}
\def \bb {\bm{b}}
\def \bz {\bm{r}}

\def \tM {\bar{M}}

\def \tA {A}
\def \tQ {Q}
\def \tR {R}
\def \tY {Y}
\def \tW {W}
\def \tZ {Z}
\def \tX {X}

\begin{document}
\nolinenumbers
\title{A fast offline/online forward solver for stationary transport equation with multiple inflow boundary conditions and varying coefficients}
\author{Jingyi Fu\thanks{School of Mathematics, Institute of Natural Sciences, Shanghai Jiao Tong University, 200240, Shanghai. Email : nbfufu@sjtu.edu.cn}
	\and
	Min Tang\thanks{School of Mathematics, Institute of Natural Sciences and MOE-LSC, Shanghai Jiao Tong University, 200240, Shanghai. 
		Email : tangmin@sjtu.edu.cn. }
}

	\maketitle
	
	\begin{abstract}
		
		It is of great interest to solve the inverse problem of stationary radiative transport equation (RTE) in optical tomography. The standard way is to formulate the inverse problem into an optimization problem, but the bottleneck is that one has to solve the forward problem repeatedly, which is time-consuming. Due to the optical property of biological tissue, in real applications, optical thin and thick regions coexist and are adjacent to each other, and the geometry can be complex. To use coarse meshes and save the computational cost, the forward solver has to be asymptotic preserving across the interface (APAL). 
		In this paper, we propose an offline/online solver for RTE. The cost at the offline stage is comparable to classical methods, while the cost at the online stage is much lower. Two cases are considered. One is to solve the RTE with fixed scattering and absorption cross sections while the boundary conditions vary; the other is when cross sections vary in a small domain and the boundary conditions change many times. The solver can be decomposed into offline/online stages in these two cases.  One only needs to calculate the offline stage once and update the online stage when the parameters vary. Our proposed solver is much cheaper when one needs to solve RTE with multiple right-hand sides or when the cross sections vary in a small domain, thus can accelerate the speed of solving inverse RTE problems. We illustrate the online/offline decomposition based on the Tailored Finite Point Method (TFPM), which is APAL on general quadrilateral meshes. 
	\end{abstract}

	\section{Introduction}
	
	
	Optical tomography (OT) is a non-invasive functional imaging of tissue to assess
	physiological function. It can detect and characterize breast cancer or other soft tissue lesions. In OT, a narrow collimated beam is sent into biological tissues, and the light that
	propagates through the medium is collected by an array of
	detectors. The sources and measurement locations are adjusted to be able to recover the material properties \cite{CBK2008,HY2016}.

	 The propagation of light in complex media can be described by the following stationary RTE 
	\eq$\bu\cdot\nabla\psi(\bz,\bu)+\big(\sigma_{S}(\bz)+\sigma_a(\bz)\big)\psi(\bz,\bu)=\sigma_{S}(\bz)\int_{S}K(\bu,\bu')\psi(\bz,\bu')\dd \bu',\label{RTE}$
	where $\psi(\bz,\bu)$ represents the photon density at position $\bz\in \Omega \subset \mathbb{R}^3$ and traveling in direction $\bu\in S$ with $S$ being an unit ball. 
	$\sigma_{S}(\bz)$, $\sigma_a(\bz)$ represent respectively the scattering cross section and absorption cross section; $K(\bu,\bu')$ is the scattering kernel that gives the probability that a particle traveling with direction $\bu'$ being scattered to direction $\bu$. 
	Boundary conditions are
	\eq$\psi(\bz,\bu)=\psi_{\Gamma^-}(\bz,\bu),\quad (\bz,\bu)\in \Gamma^-=\{\bz\in\Gamma=\pt \Omega, \quad \bu\cdot \bm{n}_{\bz}<0\}.\label{3dbc}$
	where $\bm{n}_{\bz}$ is the outward normal vector at $\bz\in \Gamma$. In order to probe the structure of highly scattering media, OT needs to solve the inverse stationary radiative transport equation (RTE), which attracts a lot of attention in the past decade. 
	As pointed in \cite{ren2010recent}, due to recent technical development, a vast number of source-detector pairs can be obtained, and it is of great interest to solve inverse RTE with substantial data sets.

	Inverse stationary RTE has been extensively studied both analytically and numerically. The uniqueness and stability results have been analyzed in \cite{Choulli_1996,choulli1999inverse}. 
	The measured data is usually a bounded linear functional of $\psi$, which can be denoted by $\mathfrak{M}\psi$. In order to get $\sigma_S(z)$, $\sigma_a(z)$ from the measured data $M$, one has to iteratively update $\sigma_S(z)$, $\sigma_a(z)$ in such a way that the forward RTE generates $\mathfrak{M}\psi$ that match $M$ with higher and higher accuracy. More precisely, one has to minimize the following objective function
	\begin{equation}\label{optimize}
		\frac{\alpha}{2}\|\mathfrak{M}\psi-M\|^2+\frac{\beta}{2}\mathfrak{R}(\sigma_a,\sigma_S)
	\end{equation}
	subject to the constraints \eqref{RTE}. Here $\alpha$, $\beta$ are two tune parameters; $\|\mathfrak{M}\psi-M\|^2$ is to quantify the difference between the model predictions and measurements; $\mathfrak{R}(\sigma_a,\sigma_S)$ is the regularization term. In this paper, we consider only inverse boundary value problems in the sense that the measurements are all taken at the boundary. There are two ways to solve the minimization problem, one is to convert \eqref{optimize} into an unconstrained optimization problem, and the other is to solve the constrained optimization problem directly \cite{ren2010recent}. In the first approach, one often first linearizes the problem around some known background to obtain a linear inverse problem. Then Green's functions that solve adjoint RTEs are needed to obtain the constraints that $\sigma_a$, $\sigma_S$ satisfy. One has to solve the forward and adjoint RTEs many times. For the second approach, forward and minimization problems are solved all at once by introducing a Lagrange multiplier. 
	The bottleneck of inverse RTE is due to the fact that the forward and adjoint RTEs depend on both spatial and angular variables. More than $90$ percent of the computational time in inverse RTE problems is for solving forward and adjoint RTEs and it is of great interest to design fast solvers for the forward problem.  
	
	An enormous amount of literature on forward solvers of steady-state RTE can be found. Two categories of methods are used: Monte Carlo methods \cite{bhan2007condensed,lee2015hybrid,hayakawa2007coupled,ueki1998kinetic} and deterministic discretization methods.
	\cite{warsa2004krylov,adams2001discontinuous,anli1996spectral,azmy1992arbitrarily,brennan2001split,lawrence1986progress,lewis1984computational,warsa2002fully,HHE2010,YCS2016}. 
	There are many challenges in the numerical simulations of RTE. The most important ones are: 1)  $\psi(\bz,\bu)$ in \eqref{RTE} depends on five independent variables (three in space and two in direction), which is costly to solve and require ample storage space; 2) the mean free path (the average distance that a particle moves between two successive collisions), varies a lot for different materials. When the mean free path is small or large, the materials are respectively called \textit{optical thick} or \textit{optical thin}. To achieve uniform convergence order in both optical thin and thick regions, asymptotic preserving (AP) schemes have to be employed. Though many schemes in the literature are not AP, there are a lot of AP schemes as well, including finite element method \cite{adams2001discontinuous}, discontinuous Galerkin method \cite{LM1989,SH2021}, finite difference method, \cite{HanTwo,LFH2019} and finite volume methods \cite{C2014,M2013}; 3)
	When optical thick and thin materials are adjacent to each other in the computational domain, there may exhibit boundary/interface layers. Meshes should be fine enough to capture the fast-changing fluxes at boundary/interface layers. However, it is not practical to resolve all layers. Then if a scheme can guarantee that the solution is valid away from the layers by using coarse meshes, we call that the scheme is asymptotic preserving across the layers (APAL). There are not many APAL schemes on general meshes in the literature. For example, to get a uniform error estimate of an upwind Discontinuous Galerkin method in slab geometry, the authors in \cite{SH2021} balance the internal discretization error with the error introduced by the unresolved boundary layer. In \cite{ShiA}, an APAL uniform convergent finite difference scheme that is valid up to the boundary and interface layer is proposed.

	Due to the optical property of biological tissue, in real applications, optical thin and thick regions coexist and are adjacent to each other, and the geometry can be complex. According to \cite{PW1989}, for an extensive range of biological tissues and different optical wavelengths,  the scattering coefficient is typically on the order of
$10^3-10^4m^{-1}$ or more, while the absorption coefficient is usually two orders of magnitude lower. Therefore, the diffusion equation is a good approximation for the forward model in OT.  However, the diffusion approximation does not hold in low-scattering regions (e.g., CSF space and trachea), highly absorbing regions (e.g., hematoma), and in the vicinity of light sources. To get better accuracy, RTE-based DOT algorithms are required, \cite{KH1999,CBK2008} and the forward RTE solvers have to be APAL on the general mesh.

	If we look more closely at the requirements for forward solvers of the above mentioned linear reconstruction and PDE-constraint approaches, we observe that \cite{ren2010recent}
	\begin{itemize}
		\item In the linear reconstruction, one solves one forward RTE for every source and one adjoint RTE for every detector. For different sources/detectors, only the boundary conditions in the forward/adjoint RTEs are different, and all other parameters are the same. Therefore linear systems with multiple right-hand sides are solved, and solutions in the whole computational domain are needed to get the constraints for $\sigma_S$, $\sigma_a$.
		\item In the nonlinear reconstruction, the minimization problem is solved by nonlinear iteration. In each iteration, forward RTEs must be solved for every source and adjoint RTEs for every measurement. The cross sections have to be updated in each iteration. In some real clinical applications, as reviewed in \cite{CBK2008,HY2016}, since OT is non-invasive, the lesion area can be much smaller than the whole computational domain that includes all sources and detectors. People only need to recover the cross sections in some small regions of interest.
	\end{itemize}

	Therefore, forward solvers that can efficiently deal with the following two cases are particularly interested in inverse RTE problems.
	\begin{itemize}
		\item \textbf{Case I.} The cross sections $\sigma_S(\bz)$, $\sigma_a(\bz)$ remain the same, boundary conditions $\psi_{\Gamma^-}(\bz,\bu)$ are chosen from a large data set. 
		\item \textbf{Case II.} The cross sections $\sigma_S(\bz)$, $\sigma_a(\bz)$ change only inside a small region or several small regions of interest, $\psi_{\Gamma^-}(\bz,\bu)$ are chosen from a large data set.  
	\end{itemize}

	If the equation has to be solved many times with different boundary conditions, one needs to design schemes that can deal with multiple right-hand sides. For example, the matrices in Diffuse Optical Tomography change slowly from one step to the next, and for each matrix, one has to solve a set of linear systems with multiple shifts and multiple right-hand sides. To reduce the total number of iterations over all linear systems, the authors \cite{kilmer2006recycling} proposed strategies for recycling Krylov subspace information. The Lanczos method has been used in \cite{saad1987} for solving symmetric linear systems with several right-hand sides. Similar ideas have been applied to calculate eigenvalue and eigenvector approximations of a Hermitian operator \cite{stathopoulos2010computing}.  Related problems have been discussed in \cite{kilmer2001qmr,kalantzis2013accelerating,zhang2020topology,bakhos2015fast}.
	The other approach is offline-online decomposition to deal with multiple right-hand sides. Some solvers can be divided into offline/online stages: only the online stage has to be updated for different sources/detectors, while the offline stage keeps the same. As far as the computational cost of the online stage is much lower than other schemes, the solver can be faster when multiple right-hand sides are solved. The idea of dividing the computation into offline/online stages is not new. Most algorithms computing elliptic PDEs on rough media (numerical homogenization) are divided into offline/online stages \cite{chen2019schwarz}. Local bases that consider the roughness of media are computed at the offline stage, and a global matrix on a coarse grid is assembled and solved at the online stage. The idea has also been extended to RTE problems based on Schwarz iteration in \cite{chen2019low}.

	We propose a new idea of dividing finite difference methods for RTE into offline/online stages. We will show that after the division, the forward solver is super-efficient in the two situations mentioned above, Case I and II. It is important to note that the solver may not be faster than other schemes when RTE is solved only once, but it will be very efficient when one solves inverse RTE with substantial data sets since the efficiency of the forward solver depends on the computational cost at the online stage. We illustrate the idea based on the Tailored finite point method (TFPM) proposed in \cite{ShiA,HanTwo,chen2018uniformly}, which is shown to be APAL. It has been demonstrated both analytically and numerically in \cite{ShiA} that the 1D TFPM for RTE is uniformly second-order convergent with respect to the mean free path up to the boundary and interface layers. 2D TFPM for RTE with isotropic and anisotropic scattering has been constructed in \cite{HanTwo,chen2018uniformly}, and the scheme can numerically capture the 2D boundary and interface layers with coarse meshes. Extension to general unstructured quadrilateral meshes has been done in \cite{WTF2022}. 
	Though TFPM has been developed for anisotropic scattering, we focus on the isotropic case in this paper to simplify the description. A similar idea can be easily extended to the anisotropic scattering.



	The organization of this paper is as follows. In section 2, we review TFPM in both 1D and 2D. The way of decomposing the solver into offline/online stages in both 1D and 2D are respectively given in sections 3 and 4. The computational costs at the offline and online stages for Case I and Case II are estimated. The numerical performance of the solver is presented in section 5, and we can see that the CPU time of the online stage is much lower than the preconditioned GMRES. Finally, we conclude with some discussions in section 6.

	\section{Review of the model and TFPM}
	\subsection{The Model}
	After nondimensionalization, RTE with isotropic scattering is
	\eq$
	\left\{\begin{aligned}
		& \bu\cdot\nabla\psi(\bz,\bu)+\left(\frac{\tilde\sigma_{S}(\bz)}{\varepsilon(\bz)}+\varepsilon(\bz)\tilde\sigma_{a}(\bz)\right)\psi(\bz,\bu)=\frac{\tilde\sigma_{S}(\bz)}{\varepsilon(\bz)}\frac{1}{|S|}\int_{S}\psi(\bz,\bu')\dd \bu', && \bz\in \Omega, \ \bu\in S, \\
		& \psi(\bz,\bu)=\psi_{\Gamma^-}(\bz,\bu), && (\bz,\bu)\in \Gamma^-. \\
	\end{aligned}\right.
	\label{3deq}$
	Here $\varepsilon(\bz)$ is a space dependent dimensionless parameter given by the ratio between mean free path and characteristic length \cite{LM1989}\blue{, which we call rescaled mean free path hereafter}. For example, bovine muscles have a scattering coefficient of about $1.7*10^4m^{-1}$ and an absorption coefficient of around $23m^{-1}$ \cite{PW1989}, then one can take the characteristic length to be $0.25cm$ and then $\sigma_s(\bz)=2.125*20*(0.25cm)^{-1}$, $\sigma_a(\bz)=1.15*0.05*(0.25cm)^{-1}$. In this case, one can take $\varepsilon(\bz)=1/20$ and the nondimensionalized scattering and absorption coefficients to be $\tilde\sigma_s(\bz)=2.125$, $\tilde\sigma_a(\bz)=1.15$. According to  \cite{ADAMS20041963,doi:10.1063/1.529374},
 when $\varepsilon\to 0$ and the inflow boundary condition $\psi_{\Gamma^-}$ is independent of $\bu$, the total density $\phi(\bz)=\frac{1}{|S|}\int_{S}\psi(\bz,\bu')\dd \bu'$ of \eqref{3deq} satisfies the following diffusion limit equation:
	\eq$-\nabla \cdot \left(\frac{1}{3\tilde\sigma_{S}(\bz)}\nabla \phi \right)+\tilde\sigma_a(\bz) \phi=0,\label{3ddiff}$
	with the boundary conditions same as $\psi_{\Gamma^-}$. For biological tissue like bovine muscles, since one can take $\varepsilon(\bz)=1/20$, which is very small, the diffusion approximation \eqref{3ddiff} is used as the forward model. It is important to note that the introduction of $\varepsilon$ is for the convenience of asymptotic analysis. It can take different values when the characteristic lengths change.

	When the equation is further reduced to slab geometry, the photon density $\psi(x,\mu)$ is defined on the domain $[x_{l},x_{r}]\times [-1,1]$. The equation becomes \cite{HanTwo}
	\begin{equation}
		\mu \frac{\pt}{\pt x}\psi(x,\mu)+\left(\frac{\tilde\sigma_{S}(x)}{\varepsilon(x)}+\varepsilon(x)\tilde\sigma_{a}\right)\psi(x,\mu)=\frac{\tilde\sigma_{S}(x)}{\varepsilon(x)}\frac{1}{2}\int_{-1}^{1}\psi(x,\mu')\dd\mu', \label{1deq}
	\end{equation}
	subject to the boundary conditions
	\begin{equation}
		\psi(x_{l},\mu)=\psi_{l}(\mu),\quad \mu>0; \quad \psi(x_{r},\mu)=\psi_{r}(\mu),\quad \mu<0. \label{1dbc}
	\end{equation}
	Then when $\varepsilon\to 0$, $\phi(x)=\frac{1}{2}\int_{-1}^1\psi(x,\mu')\dd\mu'$ satisfies the diffusion equation \cite{LARSEN1980249}
	\eq$-\frac{\pt}{\pt x}\left(\frac{1}{3\tilde\sigma_{S}(x)}\frac{\pt}{\pt x}\phi\right)+\tilde\sigma_{a} \phi=0.\label{1ddiff}$  
	The discrete ordinate method is one of the most standard method for velocity discretization. The integral term on the right-hand side of \eqref{1deq} is approximated by a weighted sum of $\psi_m(x)\approx\psi(\mu_m,x)$. The discrete-ordinate approximation of \eqref{1deq} reads
	\begin{equation}
		\mu_m\frac{\dd \psi_{m}}{\dd x}+\left(\frac{\tilde\sigma_S(x)}{\varepsilon(x)}+\varepsilon(x)\tilde\sigma_a(x)\right)\psi_{m}=\frac{\tilde\sigma_S(x)}{\varepsilon(x)} \sum_{k \in V}\omega_{k}\psi_{k}, \quad m \in V, \label{1deqd}
	\end{equation}
	and the boundary conditions are
	\begin{equation}
		\psi_{m}(x_{l})=\psi_{l}(\mu_{m}),\quad\mu_m>0; \quad \psi_{m}(x_{r})=\psi_{r}(\mu_{m}),\quad\mu_m<0. \label{1dbcd}
	\end{equation} Here $\{(\omega_m,\mu_m)\}_{m\in V}$ is the quadrature set, where the index set $V=\{1,2,\cdots,M\}$ with $M$ being an even number. We choose $\mu_{M+1-m}=-\mu_{m}>0$ and $\omega_{M+1-m}=\omega_m>0$ for $m=\frac{M}{2}+1,\cdots,M$. One typical choice is the Gaussian quadrature whose details can be found in the Appendix \cite{ShiA}.

	In 2D X-Y geometry, \eqref{3deq} writes \cite{HanTwo}
	\eq$c \frac{\pt}{\pt x}\psi+s \frac{\pt}{\pt y}\psi+\left(\frac{\tilde\sigma_{S}(x,y)}{\varepsilon(x,y)}+\tilde\sigma_a(x,y)\right)\psi=\frac{\tilde\sigma_{S}(x,y)}{\varepsilon(x,y)}\frac{1}{2\pi}\int_0^{2\pi}\int_{0}^{1}\psi(x,y,c',s')\dd \zeta' \dd \theta', \label{2deq}$
	where 
	\eq$c=(1-\zeta^2)^{\frac{1}{2}}\cos\theta, \quad s=(1-\zeta^2)^{\frac{1}{2}}\sin\theta.$
	For simplicity, we consider a rectangle domain $[x_l,x_r]\times[y_b,y_t]$ and the boundary conditions become
	\eq$
	\left\{\begin{aligned}
		&\psi(x_l,y,c,s)=\psi_{l}(y,c,s),\quad &&c>0;\quad &&\psi(x_r,y,c,s)=\psi_{r}(y,c,s),\quad &&c<0; \\
		&\psi(x,y_b,c,s)=\psi_{b}(x,c,s),\quad &&s>0;\quad &&\psi(x,y_t,c,s)=\psi_{t}(x,c,s),\quad &&s<0. \\
	\end{aligned}\right.\label{2dbc}$
	In 2D, the corresponding diffusion limit equation when $\varepsilon\to 0$ becomes
	\eq$-\frac{\pt}{\pt x}\left(\frac{1}{3\tilde\sigma_S(x,y)}\frac{\pt }{\pt x}\phi\right)-\frac{\pt}{\pt y}\left(\frac{1}{3\tilde\sigma_S(x,y)}\frac{\pt }{\pt y}\phi\right)+\tilde\sigma_a \phi=0.\label{2ddiff}$

	Discrete-ordinate approximation to \eqref{2deq} writes:
	\eq$c_m\pt_x \psi_m+s_m\pt_y \psi_m +\left(\frac{\tilde\sigma_{S}(x,y)}{\varepsilon(x,y)}+\varepsilon(x,y)\tilde\sigma_a(x,y)\right) \psi_m=\frac{\tilde\sigma_{S}(x,y)}{\varepsilon(x,y)}\sum_{k \in \bV}\bomega_k\psi_k, \quad m \in \bV, \label{eqad2}$
	where $\bV=\{1,2,\cdots,\tM\}$ and $\psi_m(x,y)\approx\psi(c_m,s_m,x,y)$. $\{c_m,s_m,\bomega_m\}_{m\in \bV}$ is the 2D quadrature set that satisfies $c_m^2+s_m^2<1$, $\bomega_m>0$. Details of how to choose $c_m,s_m,\bomega_m$ are given in Appendix A.2.
	
	Inflow boundary conditions for the discrete ordinate approximation are
	\eq$
	\left\{\begin{aligned}
		&\psi_m(x_l,y)=\psi_{l,m}(y),\quad &&c_m>0;\quad &&\psi_m(x_l,y)=\psi_{l,m}(y),\quad &&c_m<0; \\
		&\psi_m(x,y_b)=\psi_{b,m}(x),\quad &&s_m>0;\quad &&\psi_m(x,y_t)=\psi_{t,m}(x),\quad &&s_m<0, \\
	\end{aligned}\right.$
	where $\psi_{l,m}(y)$, $\psi_{r,m}(y)$, $\psi_{b,m}(x)$, $\psi_{t,m}(x)$ are respectively approximations to $\psi_l(c_m,s_m,y)$, $\psi_r(c_m,s_m,y)$, $\psi_b(c_m,s_m,x)$ and $\psi_t(c_m,s_m,x)$. 

	\subsection{TFPM in 1D}
	TFPM was first introduced in \cite{han2008tailoredsg,han2008tailoredHemker} to solve the Hemker problem and later was extended to more general singular perturbation problems of elliptic equations \cite{han2008tailoredHelmholtz,han2010tailored}, anisotropic diffusion problems \cite{TW2017}. TFPMs use the local property of the solution and thus can capture the boundary or interface layers with coarse meshes. It has been extended to 1D RTE in, \cite{ShiA} and we will review its construction in the next part.
	
		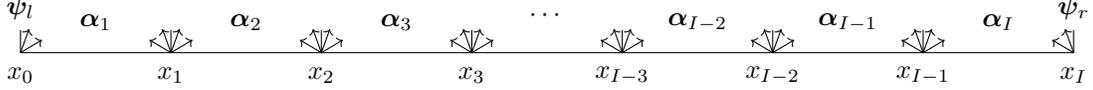
\begin{figure}[!htb]
		\centering 
		\begin{tikzpicture}
			\draw[-] (-7,0) -- (7,0);
			\foreach \i in {1,2,...,8} {\draw[-] (-9+2*\i,0)--(-9+2*\i,0.3);}
			
			\node at (-7,-0.3) {$x_0$};
			\node at (-5,-0.3) {$x_1$};
			\node at (-3,-0.3) {$x_2$};
			\node at (-1,-0.3) {$x_3$};
			\node at (1,-0.3) {$x_{I-3}$};
			\node at (3,-0.3) {$x_{I-2}$};
			\node at (5,-0.3) {$x_{I-1}$};
			\node at (7,-0.3) {$x_{I}$};
			\node at (0,0.5) {$\cdots$};
			
			\foreach \i in {1,2,...,7} {
				\draw[-{>[length=0pt 5]}] (-9+2*\i,0) -- (-9+2*\i+0.34/3,0.94/3);
				\draw[-{>[length=0pt 5]}] (-9+2*\i,0) -- (-9+2*\i+0.86/3,0.51/3);
				\draw[-{>[length=0pt 5]}] (-7+2*\i,0) -- (-7+2*\i-0.34/3,0.94/3);
				\draw[-{>[length=0pt 5]}] (-7+2*\i,0) -- (-7+2*\i-0.86/3,0.51/3);
			}
			
			\node at (-7,0.6) {$\bpsi_l$};
			\node at (7,0.6) {$\bpsi_r$};
			
			\node at (-6,0.4) {$\balpha_1$};
			\node at (-4,0.4) {$\balpha_2$};
			\node at (-2,0.4) {$\balpha_3$};
			\node at (2,0.4) {$\balpha_{I-2}$};
			\node at (4,0.4) {$\balpha_{I-1}$};
			\node at (6,0.4) {$\balpha_{I}$};
			
		\end{tikzpicture}
		\caption{Diagram of spatial and angular discretization for TFPM in 1D. Here we take $M=4$ as an example.}\label{DD1d}
	\end{figure}
	Let the grid points be $x_l=x_0<x_1<\cdots<x_I=x_r$, which include all discontinuities of functions $\tilde\sigma_a(x)$, $\tilde\sigma_T(x)=\tilde\sigma_S(x)+\varepsilon^2(x)\tilde\sigma_a(x)$, $\varepsilon(x)$.
	Then the coefficients $\sigma_a(x)$, $\sigma_T(x)$, $q(x)$, $\varepsilon(x)$ are approximated by piece-wise constants inside each cell $[x_{i-1},x_i]$ such that
	\eq$\sigma_{a,i}=\frac{\int_{x_{i-1}}^{x_i}\tilde\sigma_a(x)\dd x}{x_i-x_{i-1}},\quad \sigma_{T,i}=\frac{\int_{x_{i-1}}^{x_i}\tilde\sigma_T(x)\dd x}{x_i-x_{i-1}},\quad \varepsilon_{i}=\frac{\int_{x_{i-1}}^{x_i}\varepsilon(x)\dd x}{x_i-x_{i-1}}, \quad i=1,2,\cdots,I.
	\label{eq:psixi}$
	Then inside each cell $[x_{i-1},x_{i}]$, \eqref{1deq} can be approximated by an ODE system with constant coefficients
	\begin{equation}
		\mu_{m}\frac{d\psi_{m,i}(x)}{dx}+\frac{\sigma_{T,i}}{\varepsilon_i}\psi_{m,i}(x)=\left(\frac{\sigma_{T,i}}{\varepsilon_{i}}-\varepsilon_{i}\sigma_{a,i}\right)\sum_{k \in V}\omega_{k}\psi_{k,i}(x), \quad m \in V, \label{eqasd}
	\end{equation}
	which is equipped with boundary conditions \eqref{1dbcd} and the following interface conditions
	\eq$\psi_{m,i}(x_i)=\psi_{m,i+1}(x_i),\quad \mbox{for } m\in V, \quad i=1,2,\cdots,I-1. \label{1dicd}$
	
	The main idea of TFPM is to solve \eqref{eqasd} exactly inside each cell, and piece them together by interface conditions in \eqref{1dicd}. Let $\bpsi_{i}(x)=(\psi_{1,i}(x),\psi_{2,i}(x),\cdots,\allowbreak \psi_{M,i}(x))^T$. When $\sigma_{a,i}\neq 0$, the general solution of \eqref{eqasd} is a linear combination of the following basis functions \cite{ShiA}
	\eq$\bpsi_i^{(k)}(x)=
	\left\{\begin{aligned}
		&\bxi_i^{(k)}exp\left\{\frac{\lambda_{i}^{(k)}(x-x_{i-1})}{\varepsilon_{i}}\right\},&&\quad \lambda_{i}^{(k)}<0, \\
		&\bxi_i^{(k)}exp\left\{\frac{\lambda_{i}^{(k)}(x-x_{i})}{\varepsilon_{i}}\right\},&&\quad \lambda_{i}^{(k)}>0, \\
	\end{aligned}\right.\quad k=1,2,\cdots,M,\label{defpsik}$
	where $\lambda_o^{(k)}$ are eigenvalues of the matrix \begin{equation}
		\Lambda_i=M_{\mu}^{-1}((\sigma_{T,i}-\varepsilon^2\sigma_{a,i})W-\sigma_{T,i}I_{M}), \label{db}
	\end{equation}
	and $\bxi_{i}^{(1)},\bxi_{i}^{(2)},\cdots,\bxi_{i}^{(M)}$ are the corresponding eigenvectors. Here $M_{\mu}=diag\{\mu_{1},\mu_{2},\cdots,\mu_M\}$, $W$ is a $M \times M$ matrix with all rows being $(\omega_{1},\omega_{2},\cdots,\omega_{M})$, and $I_{M}$ is an identity matrix of size $M \times M$. Then $\bpsi^{(k)}(x)$ satisfies \eqref{eqasd} exactly. When $\sigma_{a,i}=0$, zero is a double eigenvalue of matrix $\Lambda_i$, one can let
	\eq$\bpsi_i^{(1)}(x)=\be, \qquad \bpsi_i^{(2)}(x)=\frac{\sigma_{T,i}}{\varepsilon_{i}}x\be-\bmu,$ 
	 where $\be=(1,1,\cdots,1)^T$ is a column vector of length $M$, $\bmu=(\mu_{1},\mu_{2},\cdots,\allowbreak\mu_{M})^T$. $\bpsi_i^{(1)}(x)$ and $\bpsi_i^{(2)}(x)$ satisfy \eqref{eqasd} and other basis functions are similar as in \eqref{defpsik}. In summary, general solutions to the ODE system \eqref{eqasd} write
	\begin{itemize}
		\item When $\sigma_{a,i}\neq 0$, \begin{equation}
			\bpsi_{i}(x)=
			\sum_{\lambda_i^{(k)}>0}\alpha_{k,i}\bxi_{i}^{(k)}exp\left\{\frac{\lambda_{i}^{(k)}(x-x_{i-1})}{\varepsilon_{i}}\right\}+\sum_{\lambda_i^{(k)}<0}\alpha_{k,i}\bxi_{i}^{(k)}exp\left\{\frac{\lambda_{i}^{(k)}(x-x_{i})}{\varepsilon_{i}}\right\}; \label{psin0}
		\end{equation}
		\item When $\sigma_{a,i}=0$, \begin{equation}
			\begin{aligned}
				\bpsi_{i}(x)=&\sum_{\lambda_i^{(k)}>0}\alpha_{k,i}\bxi_{i}^{(k)}exp\left\{\frac{\lambda_{i}^{(k)}(x-x_{i-1})}{\varepsilon_{i}}\right\}+\sum_{\lambda_i^{(k)}<0}\alpha_{k,i}\bxi_{i}^{(k)}exp\left\{\frac{\lambda_{i}^{(k)}(x-x_{i})}{\varepsilon_{i}}\right\} \\	&\alpha_{1,i}\be+\alpha_{2,i}\left(x\be-\frac{\varepsilon_{i}\bmu}{\sigma_{T,i}}\right). 
			\end{aligned}\label{psie0}
		\end{equation}
	\end{itemize}
	The above solution can be written into a matrix form such that
	\begin{equation}\label{eq:matrix1d}
		\bpsi_{i}(x)=A_i(x)\balpha_{i},\quad x \in [x_{i-1},x_i],
	\end{equation}
	where $\balpha_{i}=(\alpha_{1,i},\alpha_{2,i},\cdots,\alpha_{M,i})^T$ are the undetermined coefficients and $A_i(x)$ is a $M\times M$ matrices formed by fundamental solution such that
	\eq$A_i(x)=\left[\bpsi_i^{(1)} \quad \bpsi_i^{(2)} \quad \cdots\quad  \bpsi_i^{(M)}\right]. \label{defAx}$
	More details of TFPM in 1D can be found in \cite{ShiA}.

	Since \eqref{eq:matrix1d} solves \eqref{eqasd} exactly inside each cell, interface conditions in \eqref{1dicd} yield
	\begin{equation}
		\bpsi(x_i)=A_i(x_i)\balpha_{i}=A_{i+1}(x_i)\balpha_{i+1},\quad i=1,2,\cdots,I-1. \label{eqonnode}
	\end{equation} 
	Together with the $M$ equations in \eqref{1dbc} for boundary conditions, one can get the numerical solutions by solving an $IM\times IM$ linear system for $\balpha=(\balpha_1^T,\balpha_2^T,\cdots,\balpha_I^T)^T$.
	
	\subsection{TFPM in 2D}\label{sec:2.2}

	We only give the scheme construction on a rectangular domain with Cartesian grids to simplify the notations. A similar idea can be easily extended to unstructured quadrilateral meshes \cite{WTF2022}. We provide the scheme accuracy of general quadrilateral meshes in section 5.
	Let grid points in $x$ be $x_l=x_0<x_1<\cdots<x_I=x_r$, grid points in $y$ be $y_b=y_0<y_1<\cdots<y_J=y_t$.
	Cells are denoted by \eq$C_{i,j}=\{(x,y)|x_{i-1}\leq x\leq x_{i},\ y_{j-1}\leq y\leq y_{j}\},\qquad i=1,\cdots,I;j=1,\cdots,J.\label{Cij}$ 
	
	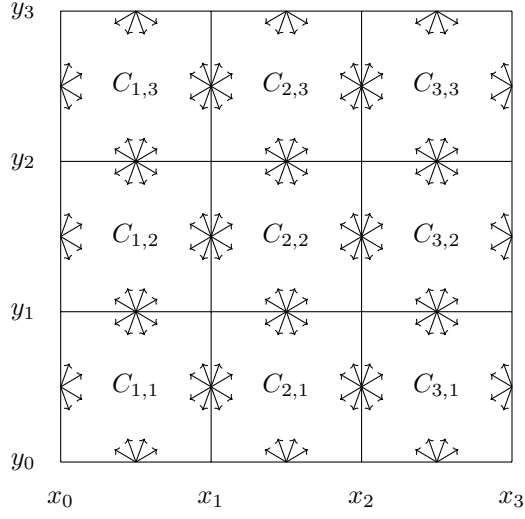
\begin{figure}[!htb]
	\centering 
	\begin{tikzpicture}
	\draw[step=2] (0,0) grid (6,6);
    \node at (1,1) {$C_{1,1}$};
    \node at (1,3) {$C_{1,2}$};
    \node at (1,5) {$C_{1,3}$};
    \node at (3,1) {$C_{2,1}$};
    \node at (3,3) {$C_{2,2}$};
    \node at (3,5) {$C_{2,3}$};
    \node at (5,1) {$C_{3,1}$};
    \node at (5,3) {$C_{3,2}$};
    \node at (5,5) {$C_{3,3}$};
    \node at (0,-0.5) {$x_0$};
    \node at (2,-0.5) {$x_1$};
    \node at (4,-0.5) {$x_2$};
    \node at (6,-0.5) {$x_3$};
    \node at (-0.5,0) {$y_0$};
    \node at (-0.5,2) {$y_1$};
    \node at (-0.5,4) {$y_2$};
    \node at (-0.5,6) {$y_3$};
	\foreach \i in {1,2} {
		\foreach \j in {1,2,3}{
			\draw[black,-{>[length=0pt 4]}] (2*\i,2*\j-1) -- (2*\i+0.34/3,2*\j-1+0.94/3);
		    \draw[black,-{>[length=0pt 4]}] (2*\i,2*\j-1) -- (2*\i+0.86/3,2*\j-1+0.51/3);
		    \draw[black,-{>[length=0pt 4]}] (2*\i,2*\j-1) -- (2*\i-0.34/3,2*\j-1+0.94/3);
		    \draw[black,-{>[length=0pt 4]}] (2*\i,2*\j-1) -- (2*\i-0.86/3,2*\j-1+0.51/3);
		    \draw[black,-{>[length=0pt 4]}] (2*\i,2*\j-1) -- (2*\i-0.34/3,2*\j-1-0.94/3);
		    \draw[black,-{>[length=0pt 4]}] (2*\i,2*\j-1) -- (2*\i-0.86/3,2*\j-1-0.51/3);
		    \draw[black,-{>[length=0pt 4]}] (2*\i,2*\j-1) -- (2*\i+0.34/3,2*\j-1-0.94/3);
		    \draw[black,-{>[length=0pt 4]}] (2*\i,2*\j-1) -- (2*\i+0.86/3,2*\j-1-0.51/3);
	    }
	}
    \foreach \i in {1,2,3} {
    	\foreach \j in {1,2}{
    	    \draw[black,-{>[length=0pt 4]}] (2*\i-1,2*\j) -- (2*\i-1+0.34/3,2*\j+0.94/3);
    	    \draw[black,-{>[length=0pt 4]}] (2*\i-1,2*\j) -- (2*\i-1+0.86/3,2*\j+0.51/3);
    	    \draw[black,-{>[length=0pt 4]}] (2*\i-1,2*\j) -- (2*\i-1-0.34/3,2*\j+0.94/3);
    	    \draw[black,-{>[length=0pt 4]}] (2*\i-1,2*\j) -- (2*\i-1-0.86/3,2*\j+0.51/3);
    	    \draw[black,-{>[length=0pt 4]}] (2*\i-1,2*\j) -- (2*\i-1-0.34/3,2*\j-0.94/3);
    	    \draw[black,-{>[length=0pt 4]}] (2*\i-1,2*\j) -- (2*\i-1-0.86/3,2*\j-0.51/3);
    	    \draw[black,-{>[length=0pt 4]}] (2*\i-1,2*\j) -- (2*\i-1+0.34/3,2*\j-0.94/3);
    	    \draw[black,-{>[length=0pt 4]}] (2*\i-1,2*\j) -- (2*\i-1+0.86/3,2*\j-0.51/3);
        }
    }
    \foreach \i in {1,2,3} {
        \draw[black,-{>[length=0pt 4]}] (2*\i-1,0) -- (2*\i-1+0.34/3,0.94/3);
        \draw[black,-{>[length=0pt 4]}] (2*\i-1,0) -- (2*\i-1+0.86/3,0.51/3);
        \draw[black,-{>[length=0pt 4]}] (2*\i-1,0) -- (2*\i-1-0.34/3,0.94/3);
        \draw[black,-{>[length=0pt 4]}] (2*\i-1,0) -- (2*\i-1-0.86/3,0.51/3);
        \draw[black,-{>[length=0pt 4]}] (2*\i-1,6) -- (2*\i-1+0.34/3,6-0.94/3);
        \draw[black,-{>[length=0pt 4]}] (2*\i-1,6) -- (2*\i-1+0.86/3,6-0.51/3);
        \draw[black,-{>[length=0pt 4]}] (2*\i-1,6) -- (2*\i-1-0.34/3,6-0.94/3);
        \draw[black,-{>[length=0pt 4]}] (2*\i-1,6) -- (2*\i-1-0.86/3,6-0.51/3);
        \draw[black,-{>[length=0pt 4]}] (0,2*\i-1) -- (0.34/3,2*\i-1+0.94/3);
        \draw[black,-{>[length=0pt 4]}] (0,2*\i-1) -- (0.86/3,2*\i-1+0.51/3);
        \draw[black,-{>[length=0pt 4]}] (0,2*\i-1) -- (0.34/3,2*\i-1-0.94/3);
        \draw[black,-{>[length=0pt 4]}] (0,2*\i-1) -- (0.86/3,2*\i-1-0.51/3);
        \draw[black,-{>[length=0pt 4]}] (6,2*\i-1) -- (6-0.34/3,2*\i-1+0.94/3);
        \draw[black,-{>[length=0pt 4]}] (6,2*\i-1) -- (6-0.86/3,2*\i-1+0.51/3);
        \draw[black,-{>[length=0pt 4]}] (6,2*\i-1) -- (6-0.34/3,2*\i-1-0.94/3);
        \draw[black,-{>[length=0pt 4]}] (6,2*\i-1) -- (6-0.86/3,2*\i-1-0.51/3);
    }
	\end{tikzpicture}\label{DD2d}
	\caption{Diagram of spatial and angular discretization for TFPM in 2D. Here we take $I=J=3$, $\tM=8$ for example.}
\end{figure}
	Assume that  $\sigma_T$, $\sigma_a$, $\varepsilon$, $q$ are continuous inside each cell,
	we approximate them by constants inside $C_{i,j}$. Then \eqref{eqad2}  can be approximated by 
	\eq$c_m\pt_x \psi_{m,i,j}+s_m\pt_y \psi_{m,i,j} +\frac{\sigma_{T,i,j}}{\varepsilon_{i,j}} \psi_{m,i,j}=\left(\frac{\sigma_{T,i,j}}{\varepsilon_{i,j}}-\varepsilon_{i,j} \sigma_{a,i,j}\right)\sum_{k \in \bV}\bomega_k\psi_{k,i,j},\quad m\in \bV,\quad (x,y)\in C_{i,j}. \label{eqasd2}$
	Solutions of different cells are pieced together by continuity of the density fluxes at cell edges, i.e. 
	$$\begin{aligned}
		&\psi_{m,i,j}(x_i,y)=\psi_{m,i+1,j}(x_i,y),\quad y\in[y_{j-1},y_{j}]; \qquad  i=1,\cdots,I-1,\ j=1,\cdots, J;\\
		&\psi_{m,i,j}(x,y_j)=\psi_{m,i,j+1}(x,y_j),\quad x\in[x_{i-1},x_{i}]; \qquad i=1,\cdots,I,\ j=1,\cdots, J-1.
	\end{aligned}
	$$


	Let $\bpsi_{i,j}(x,y)=(\psi_{1,i,j},\psi_{2,i,j},\cdots,\psi_{\tM,i,j})^T$. The TFPM proposed in \cite{HanTwo} is to approximate $\bpsi_{i,j}$ by a linear combination of basis functions and then piece solutions inside different cells together by continuity at the cell edge centers. The basis functions are of the following forms
	\begin{equation}\begin{aligned}
			&\bxi_{i,j}^{(k)}exp\left\{\frac{\lambda_{i,j}^{(k)}(x-x_{i-1})}{\varepsilon_{i,j}}\right\},\quad\mbox{for $\lambda_{i,j}^{(k)}<0$},\qquad 
			\bxi_{i,j}^{(k)}exp\left\{\frac{\lambda_{i,j}^{(k)}(x-x_{i})}{\varepsilon_{i,j}}\right\},\quad\mbox{for $\lambda_{i,j}^{(k)}>0$}\\
			&\etab_{i,j}^{(k)}exp\left\{\frac{\nu_{i,j}^{(k)}(y-y_{j-1})}{\varepsilon_{i,j}}\right\},\quad\mbox{for $\nu_{i,j}^{(k)}<0$},\qquad
			\etab_{i,j}^{(k)}exp\left\{\frac{\nu_{i,j}^{(k)}(y-y_{j})}{\varepsilon_{i,j}}\right\},\quad\mbox{for $\nu_{i,j}^{(k)}>0$.}
		\end{aligned}\label{eq:basefunction}
	\end{equation}
	Here $\lambda_{i,j}^{(k)}$
	are eigenvalues of the $\tM\times \tM$ matrix
	\begin{equation}
		\Lambda_{i,j}^c=M_{c}^{-1}((\sigma_{T,i,j}-\varepsilon_{i,j}^2\sigma_{a,i,j})W_{\tM}-\sigma_{T,i,j}I_{\tM}), \label{dbc}
	\end{equation}
	$\bxi_{i,j}^{(1)},\bxi_{i,j}^{(2)},\cdots,\bxi_{i,j}^{(\tM)}$ are corresponding eigenvectors;
	$\nu_{i,j}^{(k)}$ are eigenvalues of the $\tM\times \tM$ matrix
	\begin{equation}
		\Lambda_{i,j}^s=M_{s}^{-1}((\sigma_{T,i,j}-\varepsilon_{i,j}^2\sigma_{a,i,j})W_{\tM}-\sigma_{T,i,j}I_{\tM}, \label{dbs}
	\end{equation}
	$\etab_{i,j}^{(1)},\etab_{i,j}^{(2)},\cdots,\etab_{i,j}^{(\tM)}$ are corresponding eigenvectors. Here $M_c=diag\{c_1,\cdots,c_{\tM}\}$, $M_s=diag\{s_1,\cdots,s_{\tM}\}$,
	$W_{\tM}$ is a $\tM \times \tM$ matrix with all rows being $(\omega_{1},\omega_{2},\cdots,\omega_{\tM})$, and $I_{\tM}$ is an identity matrix of size $\tM \times \tM$.  It is easy to check that
	the basis functions in \eqref{eq:basefunction} satisfy equation \eqref{eqasd2}. When $\sigma_{a,i,j}=0$, zero is a double eigenvalue of both matrices $\Lambda_{i,j}^c$ and $\Lambda_{i,j}^s$. There are only $2\tM-1$ eigenfunctions of the form as in \eqref{eq:basefunction}. The following four basis functions that satisfy \eqref{eqasd2} are needed:
	\eq$\label{eq:sigmaa0}\begin{aligned}
		&\psi_{i,j}^{(1)}(x)=\be_{\tM}, \quad &&\psi_{i,j}^{(2)}(x,y)=\frac{\sigma_{T,i,j}}{\varepsilon_{i,j}}x\be_{\tM}-\bc, \quad \\ &\psi_{i,j}^{(3)}(x,y)=\frac{\sigma_{T,i,j}}{\varepsilon_{i,j}}y\be_{\tM}-\bs, \quad &&\psi_{i,j}^{(4)}(x,y)=\frac{\sigma_{T,i,j}}{\varepsilon_{i,j}}xy\be_{\tM}-\bs x-\bc y
		+\frac{2\varepsilon_{i,j}}{\sigma_{T,i,j}}M_c \bs,
	\end{aligned}
	$ where $\be_{\tM}$ is a $\tM\times 1$ column vector of ones, $\bc=(c_{1},c_{2},\cdots,c_{\tM})^T$, $\bs=(s_1,s_2,\cdots,s_{\tM})^T$. Therefore, when $\sigma_{a,i,j} \neq 0$, we approximate $\bpsi_{i,j}(x,y)$ by
		\begin{equation}
			\begin{aligned}
				\bpsi_{i,j}(x,y)\approx\ &\sum_{\lambda_{i,j}^{(k)}<0}\alpha_{k,i,j}\bxi_{i,j}^{(k)}exp\left\{\frac{\lambda_{i,j}^{(k)}(x-x_{i-1})}{\varepsilon_{i,j}}\right\}
				+\sum_{\lambda_{i,j}^{(k)}>0}\alpha_{k,i,j}\bxi_{i,j}^{(k)}exp\left\{\frac{\lambda_{i,j}^{(k)}(x-x_{i})}{\varepsilon_{i,j}}\right\} \\
				&+\sum_{\nu_{i,j}^{(k)}<0}\alpha_{k+\tM,i,j}\etab_{i,j}^{(k)}exp\left\{\frac{\nu_{i,j}^{(k)}(y-y_{j-1})}{\varepsilon_{i,j}}\right\}
				+\sum_{\nu_{i,j}^{(k)}>0}\alpha_{k+\tM,i,j}\etab_{i,j}^{(k)}exp\left\{\frac{\nu_{i,j}^{(k)}(y-y_{j})}{\varepsilon_{i,j}}\right\}
				,
			\end{aligned}\label{gs2d1}
		\end{equation}
		Similar approximation can be found by using \eqref{eq:sigmaa0} when $\sigma_{a,i,j}=0$ and we omit the details here. 
		
		
	Different from the 1D case, the continuity of density fluxes can  only hold at a finite number of points at the cell edges. Let $x_{i-\frac{1}{2}}=\frac{x_{i-1}+x_{i}}{2}$ for $i=1,2,\cdots,I$ and $y_{j-\frac{1}{2}}=\frac{y_{j-1}+y_j}{2}$ for $j=1,2,\cdots,J$. The approximations in \eqref{gs2d1} are pieced together by the following interface conditions:
	\eq$\begin{aligned}
		&\bpsi_{i,j}(x_i,y_{j-\frac{1}{2}})=\bpsi_{i+1,j}(x_i,y_{j-\frac{1}{2}}),&& i=1,\cdots,I-1; &&j=1,\cdots, J;\\
		&\bpsi_{i,j}(x_{i-\frac{1}{2}},y_j)=\bpsi_{i,j+1}(x_{i-\frac{1}{2}},y_j),&& i=1,\cdots,I; &&j=1,\cdots, J-1.
	\end{aligned}\label{ic2d}
	$

	As in the 1D case, \eqref{gs2d1} can be written into the following vector form 
	\begin{equation}
		\bpsi_{i,j}(x,y)=A_{i,j}(x,y)\balpha_{i,j},\quad (x,y) \in C_{i,j},
	\end{equation}
	where $\balpha_{i,j}=(\alpha_{1,i,j},\alpha_{2,i,j},\cdots,\alpha_{2\tM,i,j})^T$ are the undetermined coefficients. Interface conditions in \eqref{ic2d} and boundary conditions 
	\eq$
	\left\{\begin{aligned}
		&\psi_m(x_l,y_{j-\frac{1}{2}})=\psi_{l,m}(y_{j-\frac{1}{2}}),\ &&c_m>0,\quad &&\psi_m(x_r,y_{j-\frac{1}{2}})=\psi_{r,m}(y_{j-\frac{1}{2}}),\ &&c_m<0, \quad &&j=1,2,\cdots,J;\\
		&\psi_m(x_{i-\frac{1}{2}},y_b)=\psi_{b,m}(x_{i-\frac{1}{2}}),\ &&s_m>0,\quad &&\psi_m(x_{i-\frac{1}{2}},y_t)=\psi_{t,m}(x_{i-\frac{1}{2}}),\ &&s_m<0, \quad &&i=1,2,\cdots,I,\\
	\end{aligned}\right.\label{bc2d}$ give $2\tM I J$ equations for all $\balpha_{i,j}$.
	More details about TFPM could be found in \cite{HanTwo}. Since the fast changes at layers have been taken into account in the basis functions, TFPM has uniform convergence order even when there exhibit boundary and interface layers. 
	\begin{remark}
	In the TFPM, the discontinuities of coefficients are included in the grid points. 
	If different materials exist inside one cell, since the governing equations are different inside different materials, one can not expect good accuracy. During the iteration process of nonlinear reconstruction, since $\sigma_a$ may not be known a prior, so are the discontinuities, one may need to refine the meshes locally. Since the TFPM has uniform convergence order for general mesh as in \cite{WTF2022}, the scheme accuracy remains the same. On the other hand, as we will see later on, the offline/online decomposition remains the same when the mesh is refined locally.
	\end{remark}
	
	\section{Fast solver in 1D.}
	This part presents an efficient way of solving the large linear system constructed in section 2.2, which is adapted to multiple right-hand sides. The main idea is to build small local systems and investigate changes in the small local system when boundary conditions and parameters $\sigma_{T}$, $\sigma_a$, $\varepsilon$ vary. The construction of small local systems can be done offline, while the changes will be computed online. Therefore, the whole process can be decomposed into offline/online stages. The storage requirements and preparation time at the offline stage and the computational cost at the online stage depend on changes in two cases mentioned in the introduction.
	
	First of all, we introduce some notations. As in section 2.1, the first and last $\frac{M}{2}$ rows of $\bpsi_i(x)$ correspond to negative and positive $\mu_m$ respectively. We denote the first and last $\frac{M}{2}$ rows of $A_i(x)$ by $A_i^t(x)$ and $A_i^b(x)$ respectively. Let $$\bpsi_l^b:=(\bpsi_l(\mu_{\frac{M}{2}+1}),\bpsi_l(\mu_{\frac{M}{2}+2}),\cdots,\bpsi_l(\mu_{M}))^T,\qquad \bpsi_r^t:=(\psi_r(\mu_{1}),\psi_r(\mu_{2}),\cdots,\psi_r(\mu_{\frac{M}{2}}))^T.$$ The inflow boundary conditions in \eqref{1dbcd} write
	\eq$A_1^b(x_0)\balpha_1=\bpsi^b_{l},\label{bcalpha1}$
	\eq$A_I^t(x_I)\balpha_I=\bpsi^t_{r}.\label{bcalpha2}$
	Moreover, we denote zero matrix with size $m\times n$ by $\bm{0}_{m,n}$, and denote zero column vector with length $m$ by $\bm{0}_{m}$.
	
	\subsection{Construction of small local systems}
	We illustrate how to construct a local system of $M$ equations for each $\balpha_i$ in this subsection. To better understand the notations, diagram of spatial and angular discretization in 1D is displayed in Figure \ref{DD1d}. When $I=1$, since $A_i^b(x_0)$ and $A_i^t(x_1)$ are $\frac{M}{2} \times M$ matrices, , \eqref{bcalpha1} provides $M$ equations for $\balpha_1$. Here the coefficient matrix for $\balpha_1$ is dense, but $M$ usually is small in real applications. For example, $M=8$ to $32$ are used for 1D case in \cite{ShiA}.
	
	When $I=2$, \eqref{bcalpha1} provide $\frac{M}{2}$ equations for $\balpha_1$. On the other hand, the continuity of $\bpsi(x_1)$ provides another $M$ constraints for $\balpha_1$ such that
	\eq$A_1(x_1)\balpha_{1}=A_{2}(x_1)\balpha_{2}.\label{eq:A12}$
	Therefore $\balpha_1$ can be determined by $\balpha_2$, and the left boundary conditions in  \eqref{bcalpha1} can be expressed by $\balpha_2$ as well. More precisely,  one can eliminate $\balpha_1$ in \eqref{bcalpha1} and \eqref{eq:A12} and find 
	$\frac{M}{2}$ equations for $\balpha_2$. As far as these $\frac{M}{2}$ equations  are found, together with \eqref{bcalpha2}, one can get a $M\times M$ linear system for $\balpha_2$.
	
	Similar idea can be extended to arbitrary $I\in\mathbb{N}$.  As in Figure \ref{DD1d}, $\bpsi_l^b$ provides $\frac{M}{2}$ equations for $\balpha_1$. Since there exists a linear one-to-one map between $\balpha_{i}$ and $\balpha_{i+1}$ for $\forall i\in\{1,2,\cdots, I-1\}$,  the $\frac{M}{2}$ equations for $\balpha_i$ can be written into $\frac{M}{2}$ equations for $\balpha_{i+1}$.  By induction,  as far as the $\frac{M}{2}$ equations for $\balpha_I$ derived from the left boundary conditions  \eqref{bcalpha1} are found, together with \eqref{bcalpha2}, one can solve $\balpha_I$.
	
	Assume that the $\frac{M}{2}$ equations for $\balpha_i$ derived from the left boundary conditions are $M_i^l\balpha_i=\bb^l_i$. Similarly, one can get $\frac{M}{2}$ equations for $\balpha_i$ from the right boundary conditions denoted by $M_i^r\balpha_i=\bb_i^r$. As far as $(M_i^l$, $\bb^l_{i})$, $(M_i^r,\bb^r_{i})$ for all cells are obtained,  
	solution can be obtained by solving \eq$M_{i}\balpha_i=
	\begin{pmatrix}
	M_i^l \\
	M_i^r 
	\end{pmatrix}
	\balpha_i=
	\begin{pmatrix}
	\bb^l_{i} \\
	\bb^r_{i}  
	\end{pmatrix}=\bb_{i}.\label{ass}$ 
	inside each cell. Here $M_i^l$, $M_i^r$ are $\frac{M}{2}\times M$ matrices and $\bb^l_i$, $\bb_i^r$ are $\frac{M}{2}\times 1$ vectors. The main difficulty is to obtain $M_i$, $\bb_i$. 
	\begin{remark}\label{remark1}
		Each $M_i$ is a $M\times M$ full matrix much smaller than the whole system. One advantage of constructing small local systems is that if only solutions in a small region are needed, one can solve small systems inside the region of interest instead of the whole domain. Depending on different applications, one can decide the minimum number of small systems to solve. It is important to note that solving the small local systems in \eqref{ass}  can be done in parallel, which can significantly boost the online stage.
	\end{remark}

	\paragraph{Determine $M_i^l$, $M_i^r$, $\bb_i^l$, $\bb_i^r$ by induction.}
	Define 
	\eq$M^l_{1}\coloneqq A_1^b(x_0), \quad \bb^l_{1}\coloneqq \bpsi^b_l.\label{Mb1dl}$ In the above mentioned procedure, it is crucial to find how to determine $M_i^l$, $\bb_i^l$ by induction. Suppose that $\balpha_i$ satisfies $M_i^l\balpha_i=\bb^l_{i}$
	with $M_i^l$ being a $\frac{M}{2}\times M$ matrix and $\bb^l_{i}$ being a $\frac{M}{2}\times 1$ vector, we need to find $M^l_{i+1}$ and $\bb^l_{i+1}$ such that $M^l_{i+1}\balpha_{i+1}=\bb^l_{i+1}$. According to \eqref{eqonnode}, the most straight forward way is to use
	\eq$
	\balpha_i=\big(A_i(x_i)\big)^{-1}A_{i+1}(x_i)\balpha_{i+1}.
	\label{eq:alphai}$
	However,  when $\varepsilon_i$ tends to 0, the $M\times M$ matrix $A_i(x_i)$ is almost singular. This is because for $\lambda_{i}^{(k)}<0$, we have  $\exp\left\{\frac{\lambda_{i}^{(k)}(x_i-x_{i-1})}{\varepsilon_{i}}\right\}\to 0$,
	which indicates that all columns of $A_i(x_i)$ with $\lambda_i^{(k)}<0$ tends to $\bm{0}_{M}$.
	Therefore, it is not applicable to find $M^l_{i+1}$ and $\bb^l_{i+1}$ by substituting \eqref{eq:alphai} into $M_i^l\balpha_i=\bb^l_{i}$.

	Let 
	$$G_{i}^l=\begin{pmatrix}
		A_i(x_i) \\
		M_i^l
		\end{pmatrix}.$$ From $M_i^l\balpha_i=\bb^l_{i}$ and \eqref{eqonnode}, $\balpha_i$ satisfies
	\begin{equation}
		G_{i}^l\balpha_i= 
		\begin{pmatrix}
		A_i(x_i) \\
		M_i^l
		\end{pmatrix}
		\balpha_i
		=\begin{pmatrix}
		A_{i+1}(x_i) \\
		\bm{0}_{\frac{M}{2},M}\\
		\end{pmatrix}\balpha_{i+1}+
		\begin{pmatrix}
		\bm{0}_{M} \\
		\bb^l_{i}\\
		\end{pmatrix}
		.\label{eqb4svd}
	\end{equation}
	Let $\mathbf{l}_{i,k}$ be a $1\times \frac{3M}{2}$ vector that belongs to the left null space of $G_{i}^l$, then by left multiplying $\mathbf{l}_{i,k}$, one can eliminate $\balpha_i$ in \eqref{eqb4svd}. Since $G_{i}^l$ is a $\frac{3M}{2}\times M$ matrix, there are at least $\frac{M}{2}$ linearly independent $\mathbf{l}_{i,k}$. Hence one can find at least $\frac{M}{2}$ equations for $\balpha_{i+1}$. To get the left null space of $G_i^l$, we can use QR decomposition. Suppose that $G_i^l$ can be decomposed into $Q_i^lR_{i}^l$ with $Q_{i}^l$ being a $\frac{3M}{2}\times\frac{3M}{2}$ orthogonal matrix and $R_{i}^l$ an upper triangular matrix. Then by left multiplying $(Q_i^l)^T$ on both sides of \eqref{eqb4svd}, we get
	\begin{equation}
		R_{i}^l\balpha_i=(Q_{i}^l)^T G_{l,i} \balpha_i
		=(Q_{i}^l)^T\begin{pmatrix}
	A_{i+1}(x_i) \\
	\bm{0}_{\frac{M}{2},M}
	\end{pmatrix}
	\balpha_{i+1}+(Q_{i}^l)^T\begin{pmatrix}
	\bm{0}_{M} \\
	\bb^l_{i}\\
	\end{pmatrix}
	.\label{eqaftersvd}
	\end{equation}
	
	Here the last $\frac{M}{2}$ rows of $R_i^l$ are all zeros, thus $R_{i}^l\balpha_i$ is a $\frac{3M}{2}\times 1$ vector with the last $\frac{M}{2}$ elements being zero. We consider only the last $\frac{M}{2}$ rows of $(Q_{i}^l)^T$. Suppose \eq$(Q_{i}^l)^T=
	\begin{pmatrix}
	* &*\\
	W_{i}^l&Z_{i}^l \\
	\end{pmatrix},\label{defWZ}$ where $W_{i}^l$, $Z_{i}^l$ are respectively matrices of size $\frac{M}{2}\times M$ and $\frac{M}{2}\times \frac{M}{2}$. Then \eqref{eqaftersvd} gives
	\begin{equation}
		W_{i}^l A_{i+1}(x_i)\balpha_{i+1}+Z_{i}^l \bb^l_{i}=0,
	\end{equation}
	and $M^l_{i+1}$, $\bb^l_{i+1}$ are determined by
	\begin{equation}
		\begin{aligned}
			M^l_{i+1}&=W_{i}^l A_{i+1}(x_i), \\
			\bb^l_{i+1}&=-Z_{i}^l\bb^l_{i}.
		\end{aligned}\label{RL1dl}
	\end{equation}
	
	Similarly, let \eq$M^r_{I}\coloneqq A_{I}^t(x_I),\quad \bb^r_{I}\coloneqq \bpsi^t_r,\label{Mb1dr}$ we can find $M_i^r$ and $\bb_i^r$ $(i=1,\cdots,I-1)$ by induction.  More precisely, assume that $\balpha_i$ satisfies $M_i^r\balpha_i=\bb^r_{i}$
	for a $\frac{M}{2}\times M$ matrix $M_i^r$ and $\frac{M}{2}\times 1$ vector $\bb^r_{i}$,  one needs to find $M^r_{i-1}$ and $\bb^r_{i-1}$ such that $M^r_{i-1}\balpha_{i-1}=\bb^r_{i-1}$.  Assume that
	$$G_{i}^r=
	\begin{pmatrix}
	A_{i}(x_{i-1}) \\
	M_i^r
	\end{pmatrix}=Q_{i}^r R_{i}^r,$$ where $Q_{i}^r$ is a $\frac{3M}{2}\times \frac{3M}{2}$ orthogonal matrix  and $R_{i}^r$ is an upper triangular matrix of size $\frac{3M}{2}\times M$. Then
	\begin{equation}
		\begin{aligned}
			M^r_{i-1}&= W_{i}^r A_{i-1}(x_{i-1}), \\
			\bb^r_{i-1}&=-Z_{i}^r \bb^r_{i} \label{RL1dr}
		\end{aligned}
	\end{equation}
	satisfy $M^r_{i-1}\balpha_{i-1}=\bb^r_{i-1}$,
	where $(W_{i}^r\ Z_{i}^r)$ is the last $\frac{M}{2}$ rows of $(Q_{i}^r)^T$ with $W_{i}^r$ being of size $\frac{M}{2}\times M$ and $Z_{i}^r$ being of size $\frac{M}{2}\times \frac{M}{2}$.
	
	\begin{remark}
	The most standard Householder transformation is employed to clarify computational cost, but other more advanced methods may apply. According to \cite{golub2013matrix}, for an $m\times n(m>n)$ matrix, Householder QR factorization requires $2m n^2-\frac{2}{3}n^3$ flops, and to get the full orthogonal matrix $4m^2n-4m n^2+\frac{4}{3}n^3$ more flops are needed. Since $G_i^l$ is a $\frac{3M}{2}\times M$ matrix, it costs $\frac{20}{3}M^3$ flops to get $(Q_i^l)^T$. 
	\end{remark}
	\bigskip
	
	\subsection{Fast solver for different cases}\label{sec:3.2}
	The above approach can be decomposed into offline/online stages, two different cases are considered. We illustrate in the subsequent part how $\bb_i^l$, $\bb_i^r$ and $M_i^l$, $M_i^r$ as in \eqref{RL1dl}, \eqref{RL1dr} change with the parameters. Only the costs of updating $M_i$ and $\bb_i$ as in \eqref{ass} are displayed.
	
	\begin{itemize}[leftmargin=\widthof{Case II.}]
		\item[Case I.] {\it Influx boundary conditions $\bpsi^b_l$, $\bpsi^t_{r}$ are chosen from a large data set, $\sigma_{a,i}$, $\sigma_{T,i}$ and $\varepsilon_i$ keep the same.}
		
		Since $M_i$ is invariant in this case, we only need to update $\bb_i$ at the online stage. From \eqref{RL1dl} we have, for $i=1,\cdots,I-1$,
		\eq$\bb^l_{i+1}=-Z_{i}^l\bb^l_{i}=Z_{i}^l Z_{i-1}^l\bb^l_{i-1}=\cdots=(-1)^{i}Z_i^l Z_{i-1}^l\cdots Z_1^l\bb_1^l=(-1)^{i}Z_i^l Z_{i-1}^l\cdots Z_1^l\bpsi_l^b,$
		and similarly the second equation in \eqref{RL1dr} gives, for $i=2,\cdots,I$,
		\eq$\bb^r_{i-1}=(-1)^{I+1-i}Z_i^r Z_{i+1}^r\cdots Z_{I}^r\bpsi_r^t.$
		
		To summarize, we have 
		\paragraph{Offline/online decomposition:} 
		\begin{itemize}
			\item Offline stage. Let $M_1^l$, $M_I^r$ be as in \eqref{Mb1dl},\eqref{Mb1dr}.
			
			Compute matrices $M_i^l$, $M_i^r$ using \eqref{RL1dl} and \eqref{RL1dr}.
			
			Compute $H^l_{i}(i=1,\cdots,I-1)$, $H^r_{i}(i=2,\cdots,I)$ defined by
			\eq$H_i^l=(-1)^{i}Z_i^l Z_{i-1}^l\cdots Z_1^l, \quad H_{i}^r=(-1)^{I+1-i}Z_{i}^r Z_{i+1}^r\cdots Z_I^r. \label{defH}$ 
			
			Compute and save PLU factorization of $P_i M_i=L_i U_i$. 
			
			Matrices stored at the offline stage are $L_i$, $U_i$, $P_i$, $H_{i}^l$, $H_{i}^r$.
			\item Online stage. Compute $\bb_i^l$, $\bb_i^r$ by
			\begin{equation}\begin{aligned}
					&\bb^l_{1}=\bpsi^b_l,\qquad \bb^l_{i}=H_{i-1}^l\bpsi_l^b,\qquad\mbox{for $i=2,\cdots, I$,}\\
					&\bb^r_{I}=\bpsi^t_r,\qquad \bb^r_{i}=H_{i+1}^r\bpsi_r^t,\qquad \mbox{for $i=1,\cdots, I-1$.}\end{aligned}
				\label{blr}\end{equation}
			Solving $\balpha_i$ by
			\eq$L_i U_i \balpha_i=P_i \bb_i.$
		\end{itemize}

		\begin{remark}
		In Case I, the offline stage costs $\frac{61}{4}IM^3$ flops, online stage costs $\frac{5}{2}IM^2$ flops. The requirement of storage space that saves information from the offline stage is $\frac{3}{2}IM^2$, which is less than the nonzero elements in the assembled big matrix for solving $\balpha$. 
		The benefits of the online stage are that 1) the vectors $\bb_i^l$, $\bb_i^r$ can be updated in parallel, as well in solving the small local systems; thus, the computational time of the online stage can be independent of the number of space grids; 2) one only needs to update $\bb_i^l$, $\bb_i^r$ when solution inside the $i^{th}$ cell is needed.
		\end{remark}

		\item[Case II.]{\it $\sigma_{T,i}$, $\sigma_{a,i}$, $\varepsilon_i$ in a small subdomain $\Omega_C\subset\Omega$ ($\Omega_C$ can be unconnected), and influx boundary conditions $\bpsi^b_l$, $\bpsi^t_{r}$ are chosen from a large data set, while $\sigma_{T,i}$, $\sigma_{a,i}$, $\varepsilon_i$ in $\Omega_C^c=\Omega\setminus\Omega_C$ keep the same. }  
        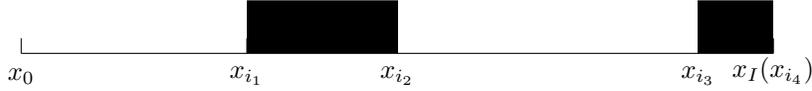
\begin{figure}[!htb]
			\centering 
			\begin{tikzpicture}
				\draw[-] (-5,0) -- (5,0);
				\foreach \i in {1,4,6,10,11} {\draw[-] (-6+\i,0)--(-6+\i,0.2);}
				\node at (-5,-0.3) {$x_0$};
				\node at (-2,-0.3) {$x_{i_1}$};
				\node at (0,-0.3) {$x_{i_2}$};
				\node at (4,-0.3) {$x_{i_3}$};
				\node at (5,-0.23) {$x_I(x_{i_4})$};
				\fill[fill=black,opacity=0.2](-2,0) rectangle (0,0.7);
				\fill[fill=black,opacity=0.2](4,0) rectangle (5,0.7);
			\end{tikzpicture}
			\caption{Example of Case II. $\sigma_{T,i}$, $\sigma_{a,i}$, $\varepsilon_i$ change in gray area but keep same in white area.}
			\label{DII}
		\end{figure}
		
		As in Figure \ref{DII}, a typical example of Case II is displayed. Here $\Omega_C=[x_{i_1},x_{i_2}]\bigcup[x_{i_3},x_{i_4}]$ and $i_4=I$. This example includes the case when $\Omega_C$ and $\Omega$ have common boundaries. 
		
		When $\sigma_{T,i}$, $\sigma_{a,i}$, $\varepsilon_i$ in $\Omega_C$ vary, by definition in \eqref{defAx}, $A_i(x_i)$, $A_i(x_{i-1})$ change. Since $M_i^l$ is determined by induction using \eqref{RL1dl}, it changes for all $i\geq i_1$. Similarly, all $M_i^r$ would change. It is expensive to update all these $M_i^l$ and $M_i^r$ and we propose a new approach in the subsequent part, in which we only update $M_i^l$, $M_i^r$ inside $\Omega_C$, while those in $\Omega_C^c$ do not change.  The idea is to construct a subsystem for all cells belong to $\Omega_C$ and the procedure can be divided into four steps:
		
		\bigskip
		\noindent \textit{First step: derive new boundary conditions for $\Omega_C$} 
		
		We use initial condition \eqref{Mb1dl} and the recursive relation \eqref{RL1dl} to obtain $M_{i_1}^l$, $\bb_{i_1}^l$. Then on node $x_{i_1}$, $M_{i_1+1}^l$, $\bb_{i_1+1}^l$ are given by
		\eq$M_{i_1+1}^l=W_{i_1}^l A_{i_1+1}(x_{i_1}), \quad \bb_{i_1+1}^l=-Z_{i_1}^l\bb_{i_1}^l=Z_{i_1}^l Z_{i_1-1}^l \bb_{i_1-1}^l=\cdots =(-1)^{i_1}Z_{i_1}^l Z_{i_1-1}^l\cdots Z_1^l \bpsi_l^b.$ Hence, we obtain $\frac{M}{2}$ equations for $\balpha_{i_1+1}$:
		\eq$W_{i_1}^l A_{i_1+1}(x_{i_1})\balpha_{i_1+1}=\hat{H}_1^l \bpsi_l^b,\label{bc1l}$
		where $\hat{H}_{1,i_1}^l=(-1)^{i_1}Z_{i_1}^l Z_{i_1-1}^l\cdots Z_1^l$. Here, $W_{i_1}^l$ and $\hat{H}_{1,i_1}^l$ are computed at the offline stage, $A_{i_1+1}(x_{i_1})$ is computed at the online stage.
		The $\frac{M}{2}$ equations of $\balpha_{i_1+1}$ in \eqref{bc1l} is regarded as new \textit{boundary condition} for $\Omega_C$ on node $x_{i_1}$. 
		
		On the other hand, node $x_I$ is a boundary point of both $\Omega_C$ and $\Omega$. New \textit{boundary condition} for $\Omega_C$ on node $x_I$ is the same as in \eqref{bcalpha2}.
		
		\bigskip
		\noindent \textit{Second step: derive new interface conditions inside $\Omega_C$} \\
		For the domain $[x_{i_2},x_{i_3}]$, the definition of $M_{i_2+j}$ $(j=1,2,\cdots,i_3-i_2)$ are different from section 3.1: we set
		\eq$M_{i_2+1}^l=A_{i_2+1}^b(x_{i_2}),\quad \bb_{i_2+1}^l=\bpsi^b(x_{i_2}),\quad M_{i_3}^r=A_{i_3}^t(x_{i_3}),\quad \bb_{i_3}^r=\bpsi^t(x_{i_3}).\label{1dicnew}$
		The relations \eqref{RL1dl} hold for $i=i_2+1,i_2+2,\cdots,i_3-1$, and \eqref{RL1dr} hold for $i=i_3,i_3-1,\cdots,i_2+2$. Then when $\sigma_a$, $\sigma_T$, $\varepsilon$ change,  $M_i^l$, $M_i^r$, $\bb_i^l$, $\bb_i^r$, $Z_i^l$, $W_i^l$, $Z_i^r$, $W_i^r$ are invariant for $i=i_2+1,i_2+2,\cdots,i_3$, since the matrices and vectors in \eqref{1dicnew} and recurrence relations \eqref{RL1dl},\eqref{RL1dr} do not change. By \eqref{RL1dl}, we can get $M_{i_3+1}^l$ and $\bb_{i_3+1}^l$ by 
		\eq$
		\begin{aligned}
		M_{i_3+1}^l=& W_{i_3}^l A_{i_3+1}(x_{i_3}), \\ \bb_{i_3+1}^l=&-Z_{i_3}^l\bb_{i_3}^l=Z_{i_3}^l Z_{i_3-1}^l\bb_{i_3-1}^l=\cdots =(-1)^{i_3-i_2}Z_{i_3}^l Z_{i_3-1}^l\cdots Z_{i_2+1}\bb_{i_2+1} \\
		&=(-1)^{i_3-i_2}Z_{i_3}^l Z_{i_3-1}^l\cdots Z_{i_2+1} \bpsi^b(i_2) =(-1)^{i_3-i_2}Z_{i_3}^l Z_{i_3-1}^l\cdots Z_{i_2+1} A_{i_2}^b(x_{i_2})\balpha_{i_2}.
		\end{aligned}$ 
		Hence, we obtain $\frac{M}{2}$ relations for $\balpha_{i_2}$ and $\balpha_{i_3+1}$:
		\eq$W_{i_3}^l A_{i_3+1}(x_{i_3})\balpha_{i_3+1}=\hat{H}_{2,i_3}^l A_{i_2}^b(x_{i_2})\balpha_{i_2},\label{ic2l}$
		with $\hat{H}_{2,i_3}^l=(-1)^{i_3-i_2}Z_{i_3}^l Z_{i_3-1}^l\cdots Z_{i_2+1}$. Here, $W_{i_3}^l$ and $\hat{H}_{2,i_3}^l$ are computed at the offline stage, $A_{i_3+1}(x_{i_3})$ and $A_{i_2}^b(x_{i_2})$ are computed at the online stage. 
		
		Similarly, we can get $M_{i_2}^r$ and $\bb_{i_2}^r$ from \eqref{RL1dr} and obtain another $\frac{M}{2}$ relations for $\balpha_{i_2}$ and $\balpha_{i_3+1}$:
		\eq$W_{i_2+1}^r A_{i_2}(x_{i_2})\balpha_{i_2}= \hat{H}_{2,i_2+1}^r A_{i_3+1}^t(x_{i_3})\balpha_{i_3+1},\label{ic2r}$
		with $\hat{H}_{2,i_2+1}^r=(-1)^{i_3-i_2}Z_{i_2+1}^r Z_{i_2+2}^r \cdots Z_{i_3}^r$. Here, $W_{i_2+1}^r$ and $\hat{H}_{2,i_2+1}^r$ are computed at the offline stage, $A_{i_2}(x_{i_2})$ and $A_{i_3+1}^t(x_{i_3})$ are computed at the online stage. 
		
		The $M$ equations of $\balpha_{i_2}$ and $\balpha_{i_3+1}$ in \eqref{ic2l} and \eqref{ic2r} are regarded as new \textit{interface conditions} for $\Omega_C$ on the two nodes $x_{i_2}$ and $x_{i_3}$.
		
		\bigskip
		\noindent \textit{Third step: solve the new system inside $\Omega_C$} \\
		Now we can build a small local system for all $\balpha_i$ inside $\Omega_C$:

		\eq$\label{3.26}\begin{pmatrix}
			W_{i_1}^l A_{i_1+1}(x_{i_1}) & & & & &\\
			& \cdots & & & &\\
			& & W_{i_2+1}^r A_{i_2}(x_{i_2}) & -\hat{H}_{2,i_2+1}^r A_{i_3+1}^t(x_{i_3}) & & \\
			& & -\hat{H}_{2,i_3}^l A_{i_2}^b(x_{i_2}) & W_{i_3}^l A_{i_3+1}(x_{i_3}) & &\\
			& & & & \cdots & \\
			& & & & &A_{I}^t(x_{I}) \\
		\end{pmatrix} \begin{pmatrix}
			\balpha_{i_1+1} \\
			\cdots \\
			\balpha_{i_2} \\
			\balpha_{i_3+1} \\
			\cdots \\ 
			\balpha_{I} 
		\end{pmatrix}=\begin{pmatrix}
			\hat{H}_{1,i_1}^l \bpsi_l^b \\
			0 \\
			0 \\ 
			0 \\
			\bpsi_r^t \\
		\end{pmatrix}$
		It has the block tri-diagonal structure, and we can apply the method described in the last section or any other solvers.
		
		\bigskip
		\noindent \textit{Final step: 
		Get the solution inside $\Omega\setminus \Omega_C$.} \\
		On $[x_0,x_{i_1}]$, we use the following boundary conditions at $x_{i_1}$
		\eq$M_{i_1}^r=A_{i_1}^t(x_{i_1}),\quad \bb_{i_1}^r=\bpsi^t(x_{i_1})=A_{i_1+1}(x_{i_1})\balpha_{i_1+1},$ and then get $M_i^r$, $\bb_i^r$ ($i< i_1$) by induction as in \eqref{RL1dr}. In such a way $M_i^r$ for $(x_{i-1},x_i)\subset[x_0,x_{i_1}]$ do not change and only the boundary conditions $\bb_i^r$ vary. Therefore, getting the solution inside $x_0,x_{i_1}$ reduces to Case I. On the other hand, after \eqref{3.26} is solved, the inflow boundary conditions of the interval $x_{i_2},x_{i_3}$ are obtained, then the solution on $[x_{i_2},x_{i_3}]$ can be found as in Case I.
		
		\bigskip
		Now we write the offline/online decomposition in general setting. Suppose $\Omega_C$ is composed of $K$ disconnected intervals $[x_{i_{2k-1}},x_{i_{2k}}]$ ($k=1,2,\cdots,K$), and there are totally $I_C(I_C\ll I)$ cells in $\Omega_C$. Here $i_k\in\{0,1,2,\cdots,I\}$, and they satisfy
		\eq$0\leq i_1<i_2<\cdots<i_{2K-1}<i_{2K}\leq I.$ 
		To simplify the notations, we suppose $i_1>0$, $i_{2K}<I$ and let $i_0=0$, $i_{2K+1}=I$. The extension to the case when $i_1=0$ or $i_{2K}=I$ is straightforward from the aforementioned discussion. Offline/online stages for Case II are:
		\paragraph{Offline/online decomposition:}
		\begin{itemize}
			\item Offline stage. For $k=0,1,\cdots,K$, in each interval $[x_{i_{2k}},x_{i_{2k+1}}]$, $\sigma_{T,i}$, $\sigma_{a,i}$, $\varepsilon_i$ do not change.
			
			Compute $M_i^l$, $M_i^r$ $(i=i_{2k}+1,i_{2k}+2,\cdots,i_{2k+1})$ by induction using \eqref{RL1dl} and \eqref{RL1dr} with
			\eq$M_{i_{2k}+1}^l=A_{i_{2k}+1}(x_{i_{2k}}),\quad M_{i_{2k+1}}^r=A_{i_{2k+1}}(x_{i_{2k+1}}).$
			
			Compute $\hat{H}_{k,i}^l$, $\hat{H}_{k,i}^r$  defined by	\eq$\hat{H}_{k,i}^l=(-1)^{i-i_{2k}}Z_{i}^l Z_{i-1}^l \cdots Z_{i_{2k}+1}^l,\quad  \hat{H}_{k,i}^r= (-1)^{i_{2k+1}+1-i}Z_{i}^r Z_{i+1}^r \cdots Z_{i_{2k+1}}^r,$
			for $i=i_{2k}+1,i_{2k}+2,\cdots,i_{2k+1}$.
			 In such setting, $G_i^l$, $G_i^r$, $Z_i^l$ and $Z_i^r$ do not change for $i=i_{2k}+1,i_{2k}+2,\cdots,i_{2k+1}$, hence $\hat{H}_{k,i}^l$ and $\hat{H}_{k,i}^r$ do not change.
			 
			 Compute and save PLU factorization of $P_i M_i=L_i U_i$ for $i=i_{2k}+1,i_{2k}+2,\cdots,i_{2k+1}$. 
			 
			Matrices stored at the offline stage are $P_i$, $L_i$, $U_i$, $\hat{H}_{k,i}^l$, $\hat{H}_{k,i}^r$, ($i=i_{2k}+1,i_{2k}+2,\cdots,i_{2k+1}$), and $W_{i_{2k+1}}^l$, $W_{i_{2k}+1}^r$ for $k=0,1,\cdots,K$.
			
			\item Online stage. 
			Compute 
			$$W_{i_1}^l A_{i_1+1}(x_{i_1}), \ \hat{H}_{1,i_1}^l \bpsi_l^b, \ W_{i_K+1}^r A_{i_K}(x_{i_K}), \ \hat{H}_{K,i_{2K}+1}^r \bpsi_r^t,$$ for the \textit{new boundary conditions} of $\Omega_C$. 
			
			For $k=2,\cdots,K$, compute 
			$$W_{i_{2k-1}}^l A_{i_{2k-1}+1}(x_{i_{2k-1}}), \ W_{i_{2k}+1}^r A_{i_{2k}}(x_{i_{2k}}), \ \hat{H}_{k,i_{2k}+1}^r A_{i_{2k-1}+1}^t(x_{i_{2k-1}}), \ \hat{H}_{k,i_{2k+1}}^l A_{i_{2k-2}}^b(x_{i_{2k-2}}),$$ 
			for the \textit{new interface conditions} at $x_{i_{2k}}$. 
			
			Solve the smaller system in $\Omega_C$.
			
			For $k=0,1,\cdots,K$, compute $\bb_i^l$, $\bb_i^r$ ($(i=i_{2k}+1,i_{2k}+2,\cdots,i_{2k+1})$) by
			\eq$
			\begin{aligned}
			& \bb_{i_{2k}+1}^l=\left\{
			\begin{aligned}
		    & \bpsi_l^b, \qquad &&  \quad \quad k=0, \\
		    & A_{i_{2k}}^r(x_{i_{2k}})\balpha_{i_{2k}}, \qquad && \quad \quad k>0, \\ \end{aligned}\right. && \bb_{i}^l=\hat{H}_{k,i-1}\bb_{i_{2k}+1}^l, && i>i_{2k}+1; \\
		    & \bb_{i_{2k+1}}^r=\left\{
			\begin{aligned}
		    & \bpsi_r^t, && k=I, \\
		    & A_{i_{2k+1}+1}^l(x_{i_{2k+1}})\balpha_{i_{2k+1}+1}, && k<I, \\ \end{aligned}\right. && \bb_{i}^r=\hat{H}_{k,i+1}\bb_{i_{2k+1}}^r, && i<i_{2k+1}, \\
			\end{aligned}
			$
			then solve $\balpha_i$ ($(i=i_{2k}+1,i_{2k}+2,\cdots,i_{2k+1})$) by $L_i U_i \balpha_i=P_i \bb_i$.
			
		\end{itemize}
		
		The new system constructed for $\balpha_i$ inside $\Omega_C$ has a similar structure as the original large sparse system. Hence methods designed for the original system can be applied to this smaller system without modification. Moreover, suppose the method used to solve the small system costs $\Upsilon(I_C,M)$ flops, the computational cost at the online stage is reduced from $\Upsilon(I,M)$ to $\frac{5}{2}(I-I_C)M^2+\Upsilon(I_C,M)$.
		
		\begin{remark}
		In Case II, the offline stage costs $\frac{61}{4} (I-I_C)M^3$ flops, online stage totally costs $\frac{5}{2}(I-I_C)M^2+\Upsilon(I_C,M)$ flops, the requirement of storage space is $\frac{3}{2}(I-I_C)M^2+\frac KM^2$.
		\end{remark}
		
	\end{itemize}
	
	\section{Fast solver in 2D}
	
	For simplicity, we only consider the rectangular domain as in section \ref{sec:2.2} and use uniform meshes, but the idea can be extended to the general case. We use the same notations as for 1D when there is no confusion.

	\subsection{Construction of small local systems using domain decomposition}
	\paragraph{Difference between 1D and 2D}
	As in 1D, we want to construct small local systems at the offline stage and compute their changes at the online stage. Recall that in 1D, the small local system is derived by induction. Using $M$ equations for $\balpha_{i-1}$ and $2M$ interface conditions of associating $\balpha_{i-1}$ and $\balpha_{i}$, one can find $M$ equations for $\balpha_i$. It is important to note that, in 1D, the dimension of $\balpha_i$ equals to the size of quadrature set, while in 2D,  the dimension of $\balpha_{i,j}$ is $2\tM$ and the size of quadrature set is $\tM$. The interface conditions between the two cells $C_{i-1,j}$ and $C_{i,j}$ are not enough to eliminate $\balpha_{i-1,j}$ from the connections between $\balpha_{i-1,j}$ and $\balpha_{i,j}$. More precisely, suppose we have $K>\tM$ equations depending only on $\balpha_{i-1,j}$, together with the $\tM$ interface conditions between the two cells $C_{i-1,j}$ and $C_{i,j}$, we have $K+\tM$ equations for $\balpha_{i-1,j}$ and $\balpha_{i,j}$. After eliminating $\balpha_{i-1,j}$, one can only find $K-\tM$ equations for $\balpha_{i,j}$, which is much less than the $K$ equations for $\balpha_{i-1,j}$. Thus the inductions as in 1D break down and we have to explore a new approach. The idea is to decompose the domain into layers and construct small local systems for each layer by induction.
	
	\paragraph{Domain decomposition}
	Suppose that the mesh is given. The first step is to divide the entire domain $\Omega$ into layers $\Omega_1$, $\Omega_2$, $\cdots$, $\Omega_S$. 
	The decomposition satisfies the following two properties:
	\begin{itemize}
		\item[(1)] each $\Omega_s$ is composed of several cells, $\Omega_1\bigcup \Omega_2\bigcup\cdots\bigcup\Omega_S=\Omega$, $\Omega_i \bigcap \Omega_j=\varnothing (1\leq i< j\leq S)$,
		\item[(2)] $\Omega_1$($\Omega_S$) has common edges only with $\Omega_2$ ($\Omega_{S-1}$), $\Omega_s$ has common cell edges only with $\Omega_{s-1}$ and $\Omega_{s+1}$ for $s=2,3,\cdots,S-1$.
	\end{itemize}
	
	In 1D, the intervals $[x_0,x_1], [x_1,x_2], \cdots, [x_{I-1},x_I]$ can be considered as a decomposition of $[x_0,x_I]$ that has the same properties.
	
	An easy way to find such decomposition is:
	\begin{enumerate}
		\item choose a subdomain $\Omega_1$ (can be unconnected), which is composed of several cells;
		\item for $s\geq 1$, let $\Omega_{s+1}$ be the set of cells that are not in $\bigcup_{s'=1}^{s}\Omega_{s'}$, and have common cell edges with at least one cell in $\bigcup_{s'=1}^{s}\Omega_{s'}$;
		\item stop with $s=S$ when $\bigcup_{s=1}^S\Omega_s=\Omega$.
	\end{enumerate}
	
\begin{figure}[!htbp]
		\centering
		\subfloat[]{\begin{tikzpicture}[scale=0.7]

				\draw[->] (-3,-3.5) -- (3,-3.5) node [anchor=west]{i};
				\draw[->] (-3.5,-3) -- (-3.5,3) node [anchor=south]{j};
				\foreach \i in {1,2,...,6} {\node at (-3.5+\i,-3.25) {\i};};
				\foreach \j in {1,2,...,6} {\node at (-3.25,-3.5+\j) {\j};};
				\draw[step=1] (-3,-3) grid (3,3);
				
				\fill[fill=red,opacity=0.7, even odd rule](-3,-3) -- (3,-3) -- (3,3) -- (-3,3) -- cycle (-2,-2) -- (2,-2) -- (2,2) -- (-2,2) -- cycle  ;
				\fill[fill=yellow!50!red,opacity=0.7, even odd rule](-2,-2) -- (2,-2) -- (2,2) -- (-2,2) -- cycle (-1,-1) -- (1,-1) -- (1,1) -- (-1,1) -- cycle;
				\fill[fill=yellow,opacity=0.7](-1,-1) rectangle  (1,1);
				\fill[fill=red,opacity=0.7] (3.5,2.25) rectangle (4,2.75); 
				\node at (4.5,2.5) {$\Omega_1$};
				\fill[fill=yellow!50!red,opacity=0.7] (3.5,-0.25) rectangle (4,0.25); 
				\node at (4.5,0) {$\Omega_2$};
				\fill[fill=yellow,opacity=0.7] (3.5,-2.75) rectangle (4,-2.25); 
				\node at (4.5,-2.5) {$\Omega_3$};
			\end{tikzpicture}\label{mesh1}}
		\hspace{1cm}
		\subfloat[]{\begin{tikzpicture}[scale=0.7]
				\draw[->] (-3,-3.5) -- (3,-3.5) node [anchor=west]{i};
				\draw[->] (-3.5,-3) -- (-3.5,3) node [anchor=south]{j};
				\foreach \i in {1,2,...,6} {\node at (-3.5+\i,-3.25) {\i};};
				\foreach \j in {1,2,...,6} {\node at (-3.25,-3.5+\j) {\j};};
				\draw[step=1] (-3,-3) grid (3,3);
				\fill[fill=red,opacity=0.7] (-2,-2) rectangle (-1,-1); 
				\fill[fill=red,opacity=0.7] (2,2) rectangle (1,1); 
				\fill[fill=yellow!50!red,opacity=0.7] (-3,-2) rectangle (-2,-1);  
				\fill[fill=yellow!50!red,opacity=0.7] (-2,-3) rectangle (-1,-2);  
				\fill[fill=yellow!50!red,opacity=0.7] (-2,-1) rectangle (-1,0);  
				\fill[fill=yellow!50!red,opacity=0.7] (-1,-2) rectangle (0,-1);  
				\fill[fill=yellow!50!red,opacity=0.7] (3,2) rectangle (2,1);  
				\fill[fill=yellow!50!red,opacity=0.7] (2,3) rectangle (1,2);  
				\fill[fill=yellow!50!red,opacity=0.7] (2,1) rectangle (1,0);  
				\fill[fill=yellow!50!red,opacity=0.7] (1,2) rectangle (0,1);  
				\fill[fill=yellow,opacity=0.7] (-3,-3) rectangle (-2,-2);
				\fill[fill=yellow,opacity=0.7] (-3,-1) rectangle (-2,0);
				\fill[fill=yellow,opacity=0.7] (-2,0) rectangle (-1,1);
				\fill[fill=yellow,opacity=0.7] (-1,-1) rectangle (0,0);
				\fill[fill=yellow,opacity=0.7] (0,-2) rectangle (1,-1);
				\fill[fill=yellow,opacity=0.7] (-1,-3) rectangle (0,-2);
				\fill[fill=yellow,opacity=0.7] (3,3) rectangle (2,2);
				\fill[fill=yellow,opacity=0.7] (3,1) rectangle (2,0);
				\fill[fill=yellow,opacity=0.7] (2,0) rectangle (1,-1);
				\fill[fill=yellow,opacity=0.7] (1,1) rectangle (0,0);
				\fill[fill=yellow,opacity=0.7] (0,2) rectangle (-1,1);
				\fill[fill=yellow,opacity=0.7] (1,3) rectangle (0,2);
				\fill[fill=green,opacity=0.7] (-3,0) rectangle (-2,1); 
				\fill[fill=green,opacity=0.7] (-2,1) rectangle (-1,2);
				\fill[fill=green,opacity=0.7] (-1,2) rectangle (0,3);
				\fill[fill=green,opacity=0.7] (-1,0) rectangle (0,1);   
				\fill[fill=green,opacity=0.7] (0,-3) rectangle (1,-2); 
				\fill[fill=green,opacity=0.7] (1,-2) rectangle (2,-1);
				\fill[fill=green,opacity=0.7] (2,-1) rectangle (3,0);
				\fill[fill=green,opacity=0.7] (0,-1) rectangle (1,0);   
				\fill[fill=blue,opacity=0.7] (-3,1) rectangle (-2,2);
				\fill[fill=blue,opacity=0.7] (-2,2) rectangle (-1,3);
				\fill[fill=blue,opacity=0.7] (1,-3) rectangle (2,-2);
				\fill[fill=blue,opacity=0.7] (2,-2) rectangle (3,-1);
				\fill[fill=black,opacity=0.7] (-3,2) rectangle (-2,3);
				\fill[fill=black,opacity=0.7] (2,-3) rectangle (3,-2);
				\fill[fill=red,opacity=0.7] (3.5,2.25) rectangle (4,2.75); 
				\node at (4.5,2.5) {$\Omega_1$};
				\fill[fill=yellow!50!red,opacity=0.7] (3.5,1.25) rectangle (4,1.75); 
				\node at (4.5,1.5) {$\Omega_2$};
				\fill[fill=yellow,opacity=0.7] (3.5,0.25) rectangle (4,0.75); 
				\node at (4.5,0.5) {$\Omega_3$};
				\fill[fill=green,opacity=0.7] (3.5,-0.75) rectangle (4,-0.25); 
				\node at (4.5,-0.5) {$\Omega_4$};
				\fill[fill=blue,opacity=0.7] (3.5,-1.75) rectangle (4,-1.25); 
				\node at (4.5,-1.5) {$\Omega_5$};
				\fill[fill=black,opacity=0.7] (3.5,-2.75) rectangle (4,-2.25); 
				\node at (4.5,-2.5) {$\Omega_6$};
			\end{tikzpicture}\label{mesh2}}
		\caption{Two examples of domain decomposition when $I=J=6$. In Figure \ref{mesh1}, $\Omega_1$ is the outermost layer of $\Omega$; in Figure \ref{mesh2}, $\Omega_1$ is the two cells at $(2,2)$ and $(5,5)$.}
		\label{meshtotal1}
	\end{figure}
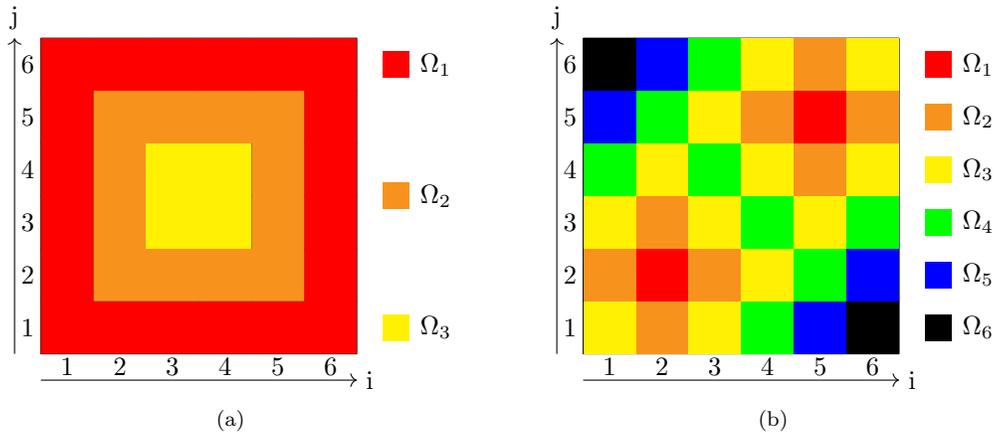

	For example, if we take $\Omega_1$ to be the set of cells on the boundary of $\Omega$, then each $\Omega_{s}$ is the $s$th outermost layer of $\Omega$, as shown in Figure \ref{mesh1}. In this type of decomposition (called annular decomposition hereafter), only cells inside $\Omega_1$ share the same cell edges with the boundary of $\Omega$.

	
	In Case II, cross sections may change in several small regions. Suppose these small regions are composed of some cells in the mesh. We can let $\Omega_1$ be these cells and determine $\{\Omega_s\}_{s=1,2,\cdots,S}$ by the above procedure. An example is shown in Figure \ref{mesh2} when $I=J=6$, $\Omega_1$ is composed of two cells at $(2,2)$ and $(5,5)$.

	 In the following part, we illustrate the way of constructing relatively small local systems for annular decomposition and give the computational cost of the offline/online stage. Then we extend it to arbitrary domain decomposition.
	
	\begin{remark}
	The computational costs vary for different domain decompositions. We choose annular decomposition for Case I due to its simplicity in description.
	\end{remark}
	
	\paragraph{Construction of small local systems}
	Suppose $I=J$ and $I$ is even. In annular decomposition, $\Omega$ can be decomposed into layers $\Omega_1$, $\Omega_2$, $\cdots$, $\Omega_{\frac{I}{2}}$ with each $\Omega_{s}$ being the $s$th outermost layer of $\Omega$ as in Figure \ref{mesh1}. The $\frac{I}{2}$th layer has 4 cells. Let $\bbeta_s$ be the column vector formed by $\balpha_{i,j}$ for all $C_{i,j}$ belonging to $\Omega_s$, then the length of $\bbeta_s$ is $8(I+1-2s)\tM$. 
	Each cell edge at the boundary gives $\frac{\tM}{2}$ equations of $\bbeta_1$; each common cell edge between two adjacent cells inside $\Omega_1$ gives $\tM$ equations of $\bbeta_1$; each common cell edge between a cell inside $\Omega_1$ and a cell inside $\Omega_2$ gives $\tM$ equations of $\bbeta_1$ and $\bbeta_2$. Hence $\bbeta_1$ satisfies
	\begin{itemize}
		\item $2I\tM$ equations from influx boundary conditions of $\Omega$, denoted by $B_1 \bbeta_1=\bb_1^B$;
		\item $4(I-1)\tM $ equations from the interface conditions between cells inside $\Omega_1$, denoted by $D_1 \bbeta_1=0$;
		\item $4(I-2)\tM$ equations from the interface conditions between $\Omega_1$ and $\Omega_2$, denoted by $\tA_1^+\bbeta_1=\tA_{2}^-\bbeta_{2}$.
	\end{itemize}

	Thus there are $(10I-12)\tM$ equations of $\bbeta_1$, $\bbeta_2$, and the size of $\bbeta_1$ is $8(I-1)\tM$. Then by eliminating $\bbeta_1$, we obtain $2(I-2)\tM$ equations of $\bbeta_2$. We next show that $2(I-2(s-1))\tM$ equations for $\bbeta_s$ can be obtained starting from boundary condition by induction. Assume that $2(I-2(s-1))\tM$ equations for $\bbeta_{s}$ have been derived from the boundary conditions by induction, denoted by $M_s^-\bbeta_s=\bb_s^-$.
	Moreover, $\bbeta_{s}$ satisfies
	\begin{itemize}
		\item $4(I+1-2s)\tM$ equations from the interface conditions between cells inside $\Omega_s$, denoted by $D_s\bbeta_s=0$; 
		\item $4(I-2s)\tM$ equations from the interface conditions of neighbouring cells between $\Omega_s$ and $\Omega_{s+1}$, denoted by $\tA_s^+\bbeta_s=\tA_{s+1}^-\bbeta_{s+1}$.
	\end{itemize}
	Hence we have
	\eq$G_s^-\bbeta_s\equiv\begin{pmatrix}
		D_s \\
		\tA_s^+\\
		M_s^- \end{pmatrix}\bbeta_{s}=\begin{pmatrix}
		\bm{0}_{4(I+1-2s)\tM,8(I-1-2s)}\\
		\tA_{s+1}^-\\
		\bm{0}_{2(I-2(s-1))\tM,8(I-1-2s)}
	\end{pmatrix}\bbeta_{s+1}+\begin{pmatrix}
		\bm{0}_{4(I+1-2s)\tM} \\
		\bm{0}_{4(I-2s)\tM}\\
		\bb_s^-
	\end{pmatrix}.$
	Here $G_s^-$ is a $2(5I+4-10s)\tM\times 8(I+1-2s)\tM$ matrix. Let the QR decomposition of $G_s^-$ be $\tQ_s^- \tR_s^-$, and $( \tY_s^- \  \tW_s^-\  \tZ_s^-)$ be the last $2(I-2s)\tM$ rows of $(\tQ_s^-)^T$, with $\tY_s^-$, $\tW_s^-$, $\tZ_s^-$ being respectively matrices of $4(I+1-2s)\tM$, $4(I-2s)\tM$ and $2(I-2(s-1))\tM$ columns. Since the last $2(I-2s)\tM$ rows of $\tR_s^-=(\tQ_s^-)^T G_s^-$ are all zeros, we get $2(I-2s)\tM$ equations for $\bbeta_{s+1}$:
	\eq$W_s^-\tA_{s+1}^-\bbeta_{s+1}+Z_s^- \bb_s^-=0.$
	Then $M_{s+1}^-$ and $\bb_{s+1}^-$ satisfying $M_{s+1}^- \bbeta_{s+1}=\bb_{s+1}^-$ can be determined by 
	\eq$\begin{aligned}
		&M_{s+1}^-=\tW_s^-\tA_{s+1}^-, \\
		&\bb_{s+1}^-=-\tZ_s^- \bb_s^-.
	\end{aligned}\label{RL2dl_annular}$
	
	On the other hand, $\Omega_{\frac{I}{2}}$ has 4 cells and $\bbeta_{\frac{I}{2}}$ satisfies
	\begin{itemize}
		\item $4\tM$ equations from the interface conditions between cells inside $\Omega_{\frac{I}{2}}$, denoted by $D_{\frac{I}{2}}\bbeta_{\frac{I}{2}}=0$;
		\item $8\tM$ equations from the interface conditions between $\Omega_{\frac{I}{2}}$ and $\Omega_{\frac{I}{2}-1}$, denoted by $\tA_{\frac{I}{2}}^-\bbeta_{\frac{I}{2}}=\tA_{\frac{I}{2}-1}^+\bbeta_{\frac{I}{2}-1}$.
	\end{itemize}
	Thus there are $12\tM$ equations for $\bbeta_{\frac{I}{2}}$, $\bbeta_{\frac{I}{2}-1}$ and size of $\bbeta_{\frac{I}{2}}$ is $8\tM$. Then by eliminating $\bbeta_{\frac{I}{2}}$, we obtain $4\tM$ equations of $\bbeta_{\frac{I}{2}-1}$.
	Similarly, assume that we have derived $2(I-2s)\tM$ equations for $\bbeta_s$, $M_s^+\bbeta_s=\bb_s^+$, since $\bbeta_{s}$ satisfies
	\begin{itemize}
		\item $4(I+1-2s)\tM$ equations from the interface conditions between cells inside $\Omega_{s}$, denoted by $D_s\bbeta_s=0$;
		\item $4(I+2-2s)\tM$ equations from the interface conditions between $\Omega_{s}$ and $\Omega_{s-1}$, denoted by $\tA_s^-\bbeta_s=\tA_{s-1}^+\bbeta_{s-1}$,
	\end{itemize}
	we have
	\eq$G_s^+\bbeta_s\equiv\begin{pmatrix}
		D_s \\
		\tA_s^-\\
		M_s^+ \end{pmatrix}\bbeta_s=\begin{pmatrix} \bm{0}_{4(I+1-2s)\tM,8(I-1-2s)}\\
		\tA_{s-1}^+\\
		\bm{0}_{2(I-2s)\tM,8(I-1-2s)}
	\end{pmatrix}\bbeta_{s-1}+\begin{pmatrix} 
		\bm{0}_{4(I+1-2s)\tM} \\
		\bm{0}_{4(I+2-2s)\tM}\\
		\bb_s^+ 
	\end{pmatrix}.$
	Here $G_s^+$ is a $2(5I+6-10s)\tM\times 4(I+1-2s)$ matrix. Let the QR decomposition of $G_s^+$ be $\tQ_s^+ \tR_s^+$, and $(Y_s^+ \ W_s^+ \ Z_s^+ )$ be the last $2(I+2-2s)\tM$ rows of $(Q_s^+)^T$, with $Y_s^+$, $W_s^+$, $Z_s^+$ being respectively matrices of $4(I+1-2s)\tM$, $4(I+2-2s)\tM$ and $2(I-2s)\tM$ columns. Then $M_{s-1}^-$ and $\bb_{s-1}^-$ satisfying $M_{s-1}^+\bbeta_s=\bb_{s-1}^+$ can be determined by
	\eq$\begin{aligned}
		&M_{s-1}^+=W_s^+\tA_{s-1}^+, \\
		&\bb_{s-1}^+=-Z_s^+ \bb_s^+,
	\end{aligned}\label{RL2dr_annular}$
    Since $\bb_{\frac{I}{2}-1}^+$ is  a zero vector, we have $\bb_s^+$ is a zero vector for $s=1,2,\cdots,\frac{I}{2}-1$.
	
	For $s=2,3,\cdots,\frac{I}{2}-1$, one can derive $2(I-2(s-1))\tM$ equations of $\bbeta_{s}$ from boundary conditions, and another $2(I-2s)\tM$ equations from $\Omega_{\frac{I}{2}}$. Together with $4(I+1-2s)\tM$ equations from the interface conditions between cells inside $\Omega_s$, one can get an $8(I+1-2s)\tM\times 8(I+1-2s)\tM$ local system for $\bbeta_s$ such that
	\eq$M_s\bbeta_s=\begin{pmatrix}
		D_s \\
		M_s^- \\
		M_s^+
	\end{pmatrix}\bbeta_s=\begin{pmatrix}
		\bm{0}_{4(I+1-2s)\tM} \\
		\bb_s^- \\
		\bm{0}_{2(I-2s)\tM} \\
	\end{pmatrix}=\bb_s.\label{ass2d_annular}$
	For $\bbeta_1$, $\bbeta_{\frac{I}{2}}$, they are solved by determined systems:
	\eq$M_1\bbeta_1=\begin{pmatrix}
		B_1 \\
		D_1 \\
		M_1^+
	\end{pmatrix}\bbeta_1=\begin{pmatrix}
		\bb_1^B \\
		\bm{0}_{4(I-1)\tM} \\
		\bm{0}_{2(I-2)\tM} \\
	\end{pmatrix}=\bb_1, \quad M_{\frac{I}{2}}\bbeta_{\frac{I}{2}}=\begin{pmatrix}
		D_{\frac{I}{2}} \\
		M_{\frac{I}{2}}^-
	\end{pmatrix}\bbeta_{\frac{I}{2}}=\begin{pmatrix}
			\bm{0}_{4\tM} \\
		\bb_{\frac{I}{2}}^- \\
	\end{pmatrix}=\bb_{\frac{I}{2}}. \label{ass2d1S_annular}$

	\begin{remark}
	In 1D case, induction using \eqref{RL1dl}, \eqref{RL1dr} starts from the boundary, and we use the boundary conditions as $M_1^l\balpha_1=\bb_1^l$, $M_I^r\balpha_I=\bb_I^r$. While in 2D case, every layer $\Omega_s$ may has cell edges on the boundary of $\Omega$ (depends on the domain decomposition). Hence, equations from influx boundary conditions of $\Omega$ are separated, and we say there is no $M_1^-$, $\bb_1^-$, $M_{I/2}^+$, $\bb_{I/2}^+$, or they are empty matrices/vectors. In the latter way, we can write \eqref{ass2d_annular} and \eqref{ass2d1S_annular} in a unified form, which is shown for general domain decomposition below.
	\end{remark}
	
	\begin{remark}
	As in 1D, we use the most standard Householder transformation to clarify the computational cost. According to \cite{golub2013matrix}, it costs $\frac{6784}{3}(I-2s)^3\tM^3+\mathcal{O}(I^2 \tM^3)$ flops to obtain $(\tQ_s^-)^T$ using Householder QR factorization, and the same flops to compute $(\tQ_s^+)^T$. Meanwhile, computing $(M_{s+1}^-, \bb_{s+1}^-)$/$(M_{s-1}^+, \bb_{s-1}^+)$ by \eqref{RL2dl_annular}/\eqref{RL2dr_annular} costs $64(I-2s)^3\tM^3+\mathcal{O}(I^2 \tM^3 )$ flops. Hence, it totally costs $\frac{1744}{3} I^4 \tM^3+\mathcal{O}(I^3 \tM^3)$ flops to construct the small local systems.
	\end{remark}

	\paragraph{General domain decomposition}
	The way to construct local systems for $\Omega_s$ in general domain decomposition is similar to the procedure in annular decomposition. As in Figure \ref{mesh2}, $\bbeta_{s}$ satisfies 
	\begin{itemize}
		\item influx boundary conditions for $\Omega$, denoted by $B_{s}\bbeta_{s}=\bb^B_{s}$;
		\item interface conditions between cells inside $\Omega_{s}$, denoted by $D_{s}\bbeta_{s}=0$;
		\item interface conditions of two cells that have common cell edges, one in $\Omega_s$ and the other in $\Omega_{s-1}$ ($s>1$), denoted by  $\tA_s^-\bbeta_s=\tA_{s-1}^+\bbeta_{s-1}$;
		\item interface conditions of two cells that have common cell edges, one in $\Omega_s$ and the other in $\Omega_{s+1}$ ($s<S$), denoted by  $\tA_s^+\bbeta_s=\tA_{s+1}^-\bbeta_{s+1}$.
	\end{itemize}
	We have the following proposition for the number of equations for $\bbeta_s$ and the degree of freedom of $\bbeta_s$:
	\begin{proposition}
		Suppose that $B_{s}$, $D_s$, $\tA_s^-$, $\tA_s^+$ have respectively $l^B_s$, $l_s^D$, $l^{A,-}_s$, $l^{A,+}_s$ rows ($l^{A,-}_1=l^{A,+}_S=0$), and the length of $\bbeta_s$ is $l^{\bbeta}_s$, then we have:
		\eq$l^{A,+}_{s-1}=l^{A,-}_{s},\qquad \mbox{for}\ s=2,3,\cdots,S,\label{l1}$
		\eq$l^{\bbeta}_s=l^B_s+l_s^D+\frac{1}{2}l^{A,-}_s+\frac{1}{2}l^{A,+}_s,\qquad \mbox{for}\ s=1,2,\cdots,S.\label{l2}$ \label{proposition1}
	\end{proposition}
	It is easy to check that Proposition \ref{proposition1} is satisfied for decompositions as in Figure \ref{mesh2}. For $\bbeta_1$, there are $l_1^B+l_1^D+l^{A,+}_1$ equations for $\bbeta_1$ and $\bbeta_2$. Thus by eliminating $\bbeta_1$ and noting \eqref{l2}, we obtain $l_1^B+l_1^D+l^{A,+}_1-l^{\bbeta}_1=\frac{1}{2}l^{A,+}_1=\frac{1}{2}l^{A,-}_2$ equations for $\bbeta_2$, which write $M_2^-\bbeta_2=\bb_2^-$. For $2\leq s \leq S-1$, assume that $\frac{1}{2}l^{A,-}_{s}$ equations of $\bbeta_s$ ($M_s^-\bbeta_s=\bb_s^-$) have been obtained starting from $\bbeta_1$, then we have $\frac{1}{2}l^{A,-}_{s}+l_s^B+l_s^D+l^{A,+}_0$ equations of $\bbeta_s$, $\bbeta_{s+1}$. Thus by eliminating $\bbeta_s$, one can get $\frac{1}{2}l^{A,-}_{s}+l_s^B+l_s^D+l^{A,+}_0-l^{\bbeta}_s=\frac{1}{2}l^{A,+}_s=\frac{1}{2}l^{A,-}_{s+1}$ equations for $\bbeta_{s+1}$. By induction, we obtain $\frac{1}{2}l^{A,-}_{s}$ equations of $\bbeta_s$ for all $s\in\{2,3,\cdots,S\}$, denoted by $M_s^-\bbeta_s=\bb_s^-$. Similarly, for all $s\in\{S-1,S-2,\cdots,1\}$, one can get $\frac{1}{2}l^{A,+}_{s}$ equations of $\bbeta_s$ by induction. We denote them by $M_s^+\bbeta_s=\bb_s^+$. Hence we have derived $l_s^D+\frac{1}{2}l^{A,-}_s+\frac{1}{2}l^{A,+}_s=l^{\bbeta}_s$ equations of $\bbeta_s$ for each $s$, and then $\bbeta_s$ can be solved by assembling the system
	\eq$M_s\bbeta_s=\begin{pmatrix}
		B_s \\
		D_s \\
		M_s^- \\
		M_s^+
	\end{pmatrix}\bbeta_s=\begin{pmatrix}
		\bb_s^B \\
		\bm{0}_{l_s^D} \\
		\bb_s^- \\
		\bb_s^+ \\
	\end{pmatrix}=\bb_s.\label{ass2d}$
	Here $M_1^-$, $M_S^+$ are empty matrices, $\bb_1^-$, $\bb_S^+$ are empty vectors. And the way of determining $M_s^-$, $\bb_s^-$, $M_s^+$, $\bb_s^+$ is similar as it in annular decomposition. For $s=1,2,\cdots,S$, let
	\eq$ G_s^-\equiv\begin{pmatrix}
		B_s \\
		D_s \\
		\tA_s^+\\
		M_s^- \end{pmatrix}=Q_s^-R_s^-, \quad G_s^+\equiv\begin{pmatrix}
		B_s \\
		D_s \\
		\tA_s^-\\
		M_s^+ \end{pmatrix}=Q_s^+R_s^+,$
	and $(X_s^{\pm}\ Y_s^{\pm}\ W_s^{\pm}\ Z_s^{\pm})$ be last $\frac{1}{2}l_s^{A,\pm}$ rows of $(Q_s^{\pm})^T$. Here $X_s^{\pm}$ have $l_s^B$ columns, $X_s^{\pm}$ have $l_s^D$ columns, $W_s^{\pm}$ have $l_s^{\mp}$ columns, $Z_s^{\pm}$ have $\frac{1}{2}l_s^{\pm}$ columns. Then 
	$M_{s+1}^-$ and $\bb_{s+1}^-$ are determined by 
	\eq$\begin{aligned}
		&M_{s+1}^-=\tW_s^-\tA_{s+1}^-, \\
		&\bb_{s+1}^-=-\tX_s^-\bb_s^B-\tZ_s^- \bb_s^-,
	\end{aligned}\label{RL2dl}$
	$M_{s-1}^+$ and $\bb_{s-1}^+$ are determined by
	\eq$\begin{aligned}
		&M_{s-1}^+=W_s^+\tA_{s-1}^+, \\
		&\bb_{s-1}^+=-\tX_s^+\bb_s^B-Z_s^+ \bb_s^+.
	\end{aligned}\label{RL2dr}$
	
	\subsection{Fast solver for different cases}
	As in 1D, we consider the two different cases discussed in the introduction. We choose suitable decomposition to reduce the cost at the online stage.
	\begin{itemize}[leftmargin=\widthof{Case II.}]
		\item  [Case I.] {\it Influx boundary conditions $\bpsi^t_b(x_{i-\frac{1}{2}})$, $\bpsi^b_t(x_{i-\frac{1}{2}})$, $\bpsi^l_r(y_{j-\frac{1}{2}})$, $\bpsi^r_l(y_{j-\frac{1}{2}})$ are chosen from a large data set, $\sigma_{a,i,j}$, $\sigma_{T,i,j}$ and $\varepsilon_{i,j}$ keep the same. }
		
		We take annular decomposition for example. It is not the best choice, but we focus on offline/online decomposition and explain the idea. In annular decomposition, when $\bpsi^t_b(x_{i-\frac{1}{2}})$, $\bpsi^b_t(x_{i-\frac{1}{2}})$,\\ $\bpsi^l_r(y_{j-\frac{1}{2}})$, $\bpsi^r_l(y_{j-\frac{1}{2}})$ vary, only $\bb_1^B$ changes. Then by \eqref{RL2dl_annular}, we have
		\eq$
		\begin{aligned}
		\bb_s^-=&-Z_{s-1}^-\bb_{s-1}^-=Z_{s-1}^-Z_{s-2}^-\bb_{s-2}^-=\cdots \\
		=&(-1)^{s-2}Z_{s-1}^-Z_{s-2}^-\cdots Z_2^-\bb_2^-=(-1)^{s-1}Z_{s-1}^-Z_{s-2}^-\cdots Z_2^- X_1^- \bb_1^B.
		\end{aligned}
		$
		On the other hand, $\bb_s^+$ is zero vector for $s=1,2,\cdots,\frac{I}{2}-1$ in annular decomposition. 
		
		As in 1D, we can save PLU factorization of $M_s$ by reducing the computational cost. The offline/online decomposition now becomes:
		\paragraph{Offline/online decomposition:}
		\begin{itemize}
			\item Offline stage. Compute matrices $M_s^-$ for $s=2,3,\cdots,\frac{I}{2}$ by induction using \eqref{RL2dl_annular} with $M_1^-$/$\bb_1^-$ being empty matrix/vector.
			Compute $H_{s}^-$ for $s>1$ by
			\eq$H_s^-=(-1)^{s-1}Z_{s-1}^-Z_{s-2}^-\cdots Z_2^- X_1^-. \label{defH2d}$
			 Compute PLU factorization of $M_s$: $P_sM_s=L_sU_s$. 
			 
			 Matrices stored at the offline stage are $L_s$, $U_s$, $P_s$, $D_s$,  $H_s^-$.
			\item Online stage. Compute $\bb_s^-$ by 
			\eq$\bb_s^-=H_s^-\bb_1^B,\quad \mbox{for}\ s=2,\cdots,S.\label{2dbm}$
			Solve $\bbeta_s$ by 
			\eq$L_s U_s \bbeta_s=P_s \bb_s.$
		\end{itemize}
		Since the $\bb_s^-$ in \eqref{2dbm} can be updated in parallel, it is straightforward to parallelize the online stage.
		\begin{remark}
		In Case I, offline stage totally costs $\frac{1768}{3}I^4 \tM^3$ flops, online stage totally costs $\frac{17}{3}I^3\tM^2$ flops, the requirement of storage space is $19I^3\tM^2$.
		\end{remark}

		\item [Case II.] {\it $\sigma_{T,i,j}$, $\sigma_{a,i,j}$, $\varepsilon_{i,j}$ in a small subdomain $\Omega_C\subset\Omega$ ($\Omega_C$ can be unconnected), and influx boundary conditions $\bpsi^t_b(x_{i-\frac{1}{2}})$, $\bpsi^b_t(x_{i-\frac{1}{2}})$, $\bpsi^l_r(y_{j-\frac{1}{2}})$, $\bpsi^r_l(y_{j-\frac{1}{2}})$ are chosen from a large data set, $\sigma_{T,i}$, $\sigma_{a,i}$, $\varepsilon_i$ in $\Omega\setminus\Omega_C$ keep the same.}  
		
		We only consider the situation when the number of spatial cells inside $\Omega_C$ is much smaller than $I\times J$. Take $\Omega_1=\Omega_C$, then when $\sigma_{T,i,j}$, $\sigma_{a,i,j}$, $\varepsilon_{i,j}$ in $\Omega_C$ vary, only $A_{i,j}(x)$ for cells inside $\Omega_1$ varies. The procedure for solving $\bbeta_s$ can be divided into four steps.
		\begin{enumerate}
		    \item Update the influence on $\bb_1$ from influx boundary conditions $\bpsi^t_b(x_{i-\frac{1}{2}})$, $\bpsi^b_t(x_{i-\frac{1}{2}})$, $\bpsi^l_r(y_{j-\frac{1}{2}})$, $\bpsi^r_l(y_{j-\frac{1}{2}})$.
		    
		    By \eqref{RL2dr} and the fact that $\bb_S^+$ is empty vector, we have
		    \eq$\begin{aligned} \bb_1^+=&-X_2^+\bb_2^B-Z_2^+\bb_2^+=-X_2^+\bb_2^B-Z_2^+(-X_3^+\bb_3^B-Z_3^+\bb_3^+)=\cdots \\
	=&-X_2^+\bb_2^B+\sum_{s=3}^{S}(-1)^{s-1} Z_2^+Z_3^+\cdots Z_{s-1}^+X_{s}^+\bb_{s}^B.\end{aligned}$
	Let 
	\eq$F_{1,2}^+=-X_2^+,\quad F_{1,s}^+=(-1)^{s-1} Z_2^+Z_3^+\cdots Z_{s-1}^+X_{s}^+, \quad 3,4,\cdots,S,$
	then $\bb_1$ can be expressed by 
	\eq$\bb_1=\begin{pmatrix}
		\bb_1^B \\
		\bm{0}_{l_1^D} \\
		\bb_s^+ \\
	\end{pmatrix}=\begin{pmatrix}
		\bb_1^B \\
		\bm{0}_{l_1^D} \\
		\sum_{s=2}^S F_{1,s}^+\bb_s^B \\
	\end{pmatrix}. \label{b1}$
		    \item Solve the local system $M_1 \bbeta_1=\bb_1$ inside $\Omega_1$. 
		    
		    \item Update the local system $M_s \bbeta_s=\bb_s$ for $s=2,3,\cdots,S$.
		    
		    Notice by \eqref{RL2dl}, $M_s^-(s>1)$ is not invariant when parameter in $\Omega_C$ vary. To avoid updating matrices, we replace \eqref{RL2dl} ($s=1$) with 
		    \eq$M_2^-=A_{2,in}^-,\quad  \bb_2^-=A_{1,out}^+\bbeta_1,\label{psicout}$
		    while compute $M_s^-(s>2)$ still with \eqref{RL2dl}. Here, equations $A_{2,in}^-\bbeta_2=A_{1,out}^+\bbeta_1$ are part of $A_2^-\bbeta_2=A_1^+\bbeta_1$, which represent the outgoing fluxes of $\Omega_C$. Then $M_s^-(s>1)$ is invariant when parameter in $\Omega_C$ vary.
		    
		    By \eqref{RL2dl} and \eqref{RL2dr}, $\bb_s^{\pm}(s>1)$ can be written by 
		    \eq$
		    \bb_s^+=-X_{s+1}^+\bb_{s+1}^B+\sum_{s'=s+2}^S(-1)^{s'-s}Z_{s+1}^+Z_{s+2}^+\cdots Z_{s'-1}^+X_{s'}^+\bb_{s'}^B,
		    $
		    \eq$\bb_s^-=-X_{s-1}^-\bb_{s-1}^B+\sum_{s'=2}^{s-2}(-1)^{s'-s}Z_{s-1}^-Z_{s-2}^-\cdots Z_{s'+1}^-X_{s'}^-\bb_{s'}^B+(-1)^{s}Z_{s-1}^-Z_{s-2}^-\cdots Z_{2}^-\bb_2^-.$
		
		Let 
		\eq$F_{s,s+1}^+=-X_{s+1}^+,\quad F_{s,s'}^+=(-1)^{s'-s} Z_{s+1}^+Z_{s+2}^+\cdots Z_{s'-1}^+X_{s'}^+, \quad s'=s+2,s+3,\cdots,S,\label{Fp}$
		\eq$F_{s,s-1}^-=-X_{s-1}^-,\quad F_{s,s'}^+=(-1)^{s'-s} Z_{s-1}^-Z_{s-2}^-\cdots Z_{s'+1}^-X_{s'}^-, \quad s'=2,3,\cdots,s-2,\label{Fm}$
		\eq$H_s=(-1)^{s}Z_{s-1}^-Z_{s-2}^-\cdots Z_{2}^-A_{1,out}^+,\label{H2}$
		then $\bb_s(s>1)$ can be expressed by 
		\eq$\bb_s=\begin{pmatrix}
		\bb_s^B \\
		\bm{0}_{l_1^D} \\
		\bb_s^- \\
		\bb_s^+
	\end{pmatrix}=\begin{pmatrix}
		\bb_s^B \\
		\bm{0}_{l_s^D} \\
		\sum_{s'=s+1}^S F_{s,s'}^+\bb_{s'}^B \\
		\sum_{s'=2}^{s-1} F_{s,s'}^-\bb_{s'}^B+H_s\bbeta_1
	\end{pmatrix}.\label{bs}$
	
	\item Solve the local system $M_s \bbeta_s=\bb_s$ inside $\Omega_s$ for $s=2,3,\cdots, S$.
	\end{enumerate}
		To sum up, we have 
		\paragraph{Offline/online decomposition:}
		\begin{itemize}
			\item Offline stage. Compute matrices $M_s^-$ by induction using \eqref{RL2dl} for $s=3,\cdots,S$ with $M_2^-=A_{2,in}^-$.
			
			Compute matrices $M_s^+$ by induction using \eqref{RL2dl} for $s=S-1,\cdots,2$ with $M_S^+$ being a $0\times l_S^{\bbeta}$ empty matrix.
			
			Compute $F_{s,s'}^+$ by \eqref{Fp} for $s=1,2,\cdots,S-1$.
			
			Compute $F_{s,s'}^-$ and $H_s$ by \eqref{Fm} and \eqref{H2} for $s=3,\cdots,S$.
			
			Matrices stored at the offline stage are $W_2^+$, $M_s^{\pm}$ for $s=2,3,\cdots,S$, $F_{s,s'}^+$ for $s=1,2,\cdots,S-1$, $F_{s,s'}^-$ and $H_s$ for $s=3,\cdots,S$.
			
			\item Online stage. Compute $\bb_1$ by \eqref{b1}, $M_1^+$ by $M_1^+=W_2^+A_1^+$, then solve $\bbeta_1$ by $M_1\bbeta_1=\bb_1$.
			
			Compute $\bb_s(s>1)$ by \eqref{bs}, then solve $\bbeta_s(s>1)$ by $M_s\bbeta_s=\bb_s$.
		\end{itemize}

    \begin{remark}
        Since the computational costs vary for different $\Omega_C$, we do not discuss the flops in this case. When the number of cells is much smaller than $I\times J$, the offline stage costs $\mathcal{O}(I^4 \tM^3)$ flops, costs at the online stage are $\mathcal{O}(I^3\tM^2)$ flops, the  requirement of storage space is $\mathcal{O}(I^3\tM^2)$.
    \end{remark}
	\end{itemize}

	\section{Numerical examples}
	
	In this section, several numerical examples are displayed to validate the accuracy and efficiency of our algorithm. In both 1D and 2D, we demonstrate the APAL property of TFPM with examples whose solutions exhibit boundary and interface layers. Uniform quadratic convergence on non-uniform meshes can be observed numerically, even when the boundary and interface layers coexist. The scheme efficiency is illustrated by the runtime of online and offline stages for different rescaled  mean free paths $\varepsilon$, numbers of spatial cells $M$, and number of spatial cells $I$ in Case I/II.  
	
	Computations shown below are performed single-threaded on an Inter Xeon Processor (Skylake, IBRS) @ 2.39 GHz, coded in Matlab. Restarted GMRES solver with both block-diagonal right preconditioner and ILU right preconditioner is chosen for comparison. Runtime at the offline stage ($T_{off}$), online stage ($T_{on}$), and GMRES ($T_{GMRES}$) are shown for both 1D and 2D examples whose solutions exhibit boundary and interface layers. For the sake of fairness, the time of matrices and vectors construction is not included in the runtime, while time for preconditioning is contained. Moreover, the tolerance of GMRES is $10^{-10}$, the tolerance of ILU preconditioner is $10^{-6}$, and GMRES restarts every 5 inner iterations. For fairness, we do not parallelize the code, but our online stage can be easily accelerated by parallelization.

	\subsection{1D Case}
	
\begin{figure}[!htbp]
	    \centering \subfloat[]{\includegraphics[width=0.5\linewidth]{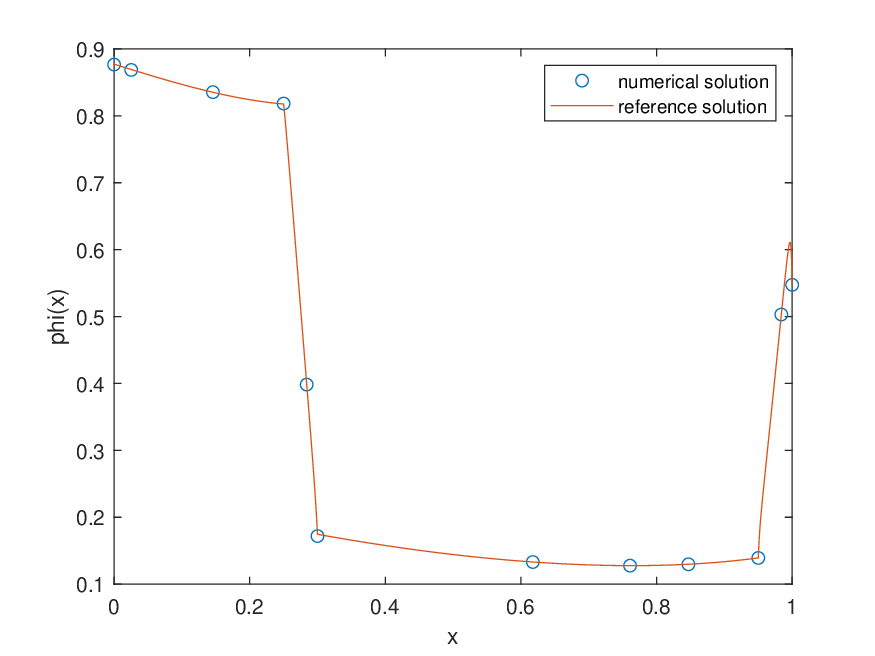}\label{1dCaseII}}  \subfloat[]{\includegraphics[width=0.5\linewidth]{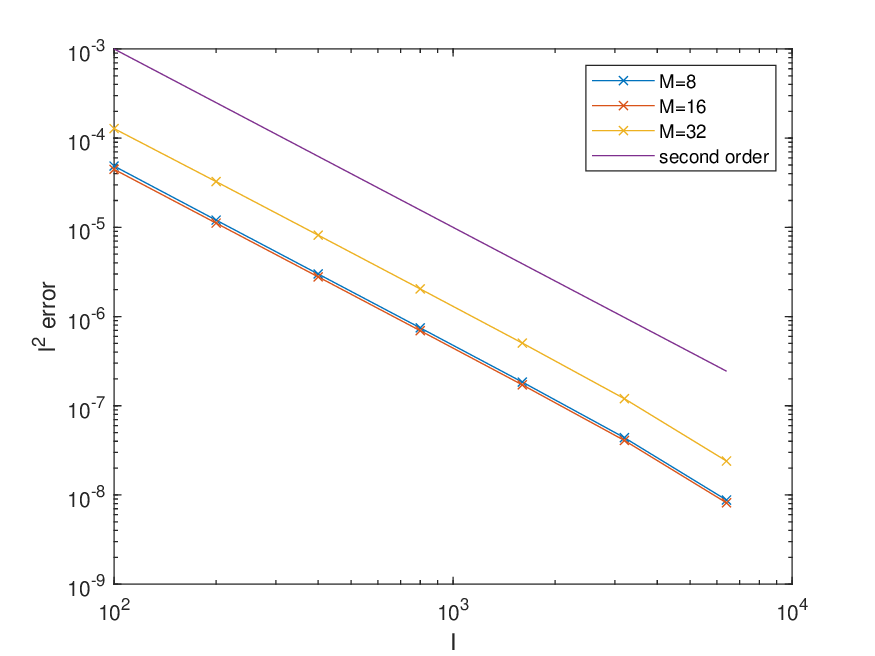}\label{err1d}}
	    \caption{Example 1. a) Density flux $\phi(x)=\sum_{k\in V}\omega_k \psi_k(x)$ with number of spatial cells $M=6$. Numerical result with non-uniform coarse mesh (circles) and the \textit{reference} solution (solid lines) are displayed. b) Convergence order with non-uniform mesh for number of discrete ordinates $M=8,16,32$, and
	    number of spatial cells $I=100,200,400,800,1600,3200,6400$.}
	    \label{figure1}
	\end{figure}

	\paragraph{Example 1: } In order to show the APAL property of TFPM in 1D, we choose $\sigma_T$, $\sigma_a$, $q$ depend on space and the magnitude of $\varepsilon$ varies. Moreover, the inflow boundary condition is chosen to be anisotropic so that the solution exhibits both boundary and interface layers. More precisely, let  
	    \begin{equation*}
	    \begin{aligned}
	    & x\in \Omega=[0,1],\quad \psi_l^b=(1,1,\cdots,1)^T,\quad \bpsi_r^t=(\mu_{M/2+1},\mu_{M/2+2},\cdots,\mu_{M})^T, \\
	    & \sigma_{T}=x^2+1, \quad \sigma_{a}=x^2+0.5, \quad \varepsilon=0.01, \quad x\in \Omega_C= [0.25,0.3] \cup [0.95,1], \\
	    & \sigma_{T}=1, \quad \sigma_{a}=x,\quad \varepsilon=1, \quad x\in \Omega\setminus \Omega_C= [0,0.25)\cup(0.3,0.95).
	    \end{aligned}
		\end{equation*}
	\paragraph{APAL property.} The interface points $0.25$, $0.3$, $0.95$ are included in the set of grid points, while other nodes of the coarsest mesh are chosen randomly from $[x_l,x_r]$. Finer meshes are refined based on the coarsest mesh. For example, the second coarsest mesh includes grid points of the coarsest mesh and the midpoint of each cell. The \textit{reference solution} $\psi_m^{exact}(x)$ refers to the result computed by the same method with uniform mesh and $\Delta x=1/12800$. As shown in Figure \ref{1dCaseII}, the solution exhibits boundary and interface layers. Even if the boundary/interface layer is not resolved numerically, the proper solution behavior can be captured with a coarse non-uniform mesh. 
	     
	The $l^2$ error between the \textit{reference solution} and numerical solution with number of spatial cells $I$ is defined by: 
	$$\sqrt{\frac{1}{IM}\sum_{i=0}^{I}\sum_{m=1}^{M}\left(\psi_{m}(x_i)-\psi_{m}^{exact}(x_i)\right)^2}.$$
	$l^2$ error for different numbers of spatial cells $I$ and numbers of discrete ordinates $M$ are shown in Figure \ref{err1d}. We can observe second-order convergence even when boundary and interface layers coexist.
	
	\paragraph{Runtime in Case I}
		We check the dependence of runtime at the offline stage ($T_{off}$), online stage ($T_{on}$) and GMRES method ($T_{gmres}$) on rescaled mean free path $\varepsilon$ inside $\Omega_C$, number of discrete ordinates $M$ and number of spatial cells $I$. In Case I, only the boundary conditions vary. Firstly, we fix $M=16$, $I=10000$. We consider different rescaled mean free paths $\varepsilon=1,0.1,0.01,0.001$ inside $\Omega_C$, $T_{off}$, $T_{on}$, and $T_{gmres}$  are shown in Table \ref{tbvarepsilon}. The results show that, the computational costs of all stages are independent of $\varepsilon$.
	
		In Table \ref{ItbM}, we fix $\varepsilon=1$, $I=10000$, and show the runtime of different $M$. $T_{off}$ increases with $M$ almost quadratically, slower than $M^3$ as in the standard Householder method. The increasing rate of $T_{on}$ is about quadratically in $M$ as well. 
		In Table \ref{ItbI}, we fix $\varepsilon=1$, $M=16$, and show the runtime of different $I$. Runtime increases linearly w.r.t $I$. We can see that the overall computational time is longer but comparable to the preconditioned GMRES method, but the online stage needs much less time. 
	    
	\begin{table}[!htbp]
	    \centering
	    \subfloat[]{\begin{tabular}{|c|c|c|c|}
			\hline
			$\varepsilon$ & $T_{off}$ & $T_{on}$ & $T_{gmres}$  \\
			\hline
			1 & 0.81 & 0.13 & 0.49 \\
		    \hline
			$10^{-1}$ & 0.82& 0.13 & 0.51\\
			\hline
			$10^{-2}$ & 0.74& 0.13 & 0.51\\
			\hline
			$10^{-3}$ & 0.75& 0.13 & 0.52\\
		    \hline
			$10^{-4}$ & 0.73& 0.13 & 0.52\\
			\hline 
			\end{tabular}		\label{tbvarepsilon}}
		\hspace{0.5cm}
		\subfloat[]{\begin{tabular}{|c|c|c|c|}
		    \hline
			$M$ & $T_{off}$ & $T_{on}$ & $T_{gmres}$\\
			\hline
			16 & 0.81 & 0.13 & 0.49 \\
			\hline
			32 & 1.91& 0.20 & 1.74\\
			\hline
			64 & 9.72 & 0.50 & 6.93\\
			\hline
			128 & 38.67 & 1.44 & 32.68\\
			\hline
			256 & 184.97 & 6.18 & 171.51\\
			\hline
		\end{tabular}
		\label{ItbM}}
	    \hspace{0.5cm}
		\subfloat[]{\begin{tabular}{|c|c|c|c|}
		    \hline
			$I$ & $T_{off}$ & $T_{on}$ & $T_{gmres}$ \\
			\hline
			10000 & 0.81 & 0.13 & 0.49 \\
			\hline
			20000 & 1.35 & 0.25& 0.98 \\
			\hline
			40000 & 2.65 & 0.51& 1.98 \\
			\hline
			80000 & 6.18 & 1.02& 4.14 \\
			\hline
			160000 & 12.05 & 1.99& 8.54 \\
			\hline
		\end{tabular}
		\label{ItbI}}
		\caption{Example 1. 1D Case I. Run time (seconds) of different stages in the fast solver (offline/online) with different rescaled mean free paths $\varepsilon$ (Table \ref{tbvarepsilon}), number of discrete ordinates $M$ (Table \ref{ItbM}) or number of spatial cells $I$ (Table \ref{ItbI}).}
		\label{Itb}
	\end{table}
	
	\paragraph{Runtime in Case II}
	In Case II, parameters $\sigma_T$, $\sigma_a$, $\varepsilon$ in $\Omega_C$
		and boundary conditions vary. Note that the size of $\Omega_C$ is one-tenth of $\Omega$ in this example. We check how the runtime depends on the number of discrete ordinates $M$ and number of spatial cells $I$.  
		 In Table \ref{IItb}, the runtime at different stages using different $M$ and $I$ are displayed. We observe that, $T_{off}$ and $T_{on}$ increase almost quadratically w.r.t $M$, and linearly w.r.t $I$. The overall computational time is larger than the preconditioned GMRES method, but the online stage needs much less time. 
		
		\begin{table}[!htbp]
		\centering
		\subfloat[]{\begin{tabular}{|c|c|c|c|}
				\hline
				$M$ & $T_{off}$ & $T_{on}$ & $T_{gmres}$\\
				\hline
				16 & 0.71 & 0.18 & 0.48 \\
				\hline
				32 & 1.71& 0.37 & 1.76 \\
				\hline
				64 & 8.54 & 1.32& 7.05\\
				\hline
				128 & 35.01 & 6.36 & 32.31 \\
				\hline
			    256 & 163.15 & 32.84 & 173.86 \\
				\hline
			\end{tabular}
			\label{IItbM}}
		\hspace{1cm}
		\subfloat[]{\begin{tabular}{|c|c|c|c|}
				\hline
				$I$ & $T_{off}$ & $T_{on}$ & $T_{gmres}$ \\
				\hline
				10000 & 0.71 & 0.18 & 0.48 \\
				\hline
				20000 & 1.34& 0.31& 0.96\\
				\hline
				40000 & 2.74& 0.64& 1.96\\
				\hline
				80000 & 4.80& 1.41& 4.12\\
				\hline
				160000 & 11.08 & 2.94 & 8.40 \\
				\hline
			\end{tabular}
			\label{IItbI}}
		\caption{Example 1. 1D Case II. Run time (seconds) of different stages in the fast solver (offline/online) with different numbers of discrete ordinates $M$ or number of spatial cells $I$.}
		\label{IItb}
	\end{table}

	\subsection{2D case}

    	\paragraph{Example 2: } We show the uniform second-order convergence of TFPM up to the boundary layer on general quadrilateral meshes. Let
    \begin{equation*}
    (x,y)\in \Omega= [0,1]\times[0,1],\quad \sigma_T(x,y)=1, \quad \sigma_a(x,y)=0.5.
    \end{equation*}
    We consider the following exact solution: 
    \eq$\psi(x,y)=\bm{\zeta} exp\left\{\frac{\frac{1}{2}\tau(x-1)+
    \frac{\sqrt{3}}{2}\tau(y-1)}{\varepsilon}\right\}+\bxi exp\left\{\frac{\lambda(x-1)}{\varepsilon}\right\}+\bm{\eta} exp\left\{\frac{\nu(y-1)}{\varepsilon}\right\}, \label{es}$
    where $\tau$ is the smallest positive eigenvalue of $(\frac{1}{2}M_c+\frac{\sqrt{3}}{2}M_s)^{-1}((\sigma_{T}-\varepsilon^2\sigma_{a})W_{\tM}-\sigma_{T}I_{\tM})$, 
    $\bm{\zeta}$ is the corresponding scaled eigenvector such that $\Vert \bm{\zeta} \Vert_{\infty}=1$; $\lambda$ is the second smallest positive eigenvalue of $M_c^{-1}((\sigma_{T}-\varepsilon^2\sigma_{a})W_{\tM}-\sigma_{T}I_{\tM})$, $\bm{\xi}$ is the corresponding scaled eigenvector such that $\Vert \bm{\xi} \Vert_{\infty}=1$; $\nu$ is the second smallest positive eigenvalue of $M_s^{-1}((\sigma_{T}-\varepsilon^2\sigma_{a})W_{\tM}-\sigma_{T}I_{\tM})$, $\bm{\eta}$ is the corresponding scaled eigenvector such that $\Vert \bm{\eta} \Vert_{\infty}=1$. The first eigenfunction is 
    not base function in the scheme construction, while the last two functions are base functions. The inflow boundary conditions are determined by the exact solution \eqref{es}.
    
\begin{figure}[!htbp]
	    \centering
	    \subfloat[]{\includegraphics[width=0.5\linewidth]{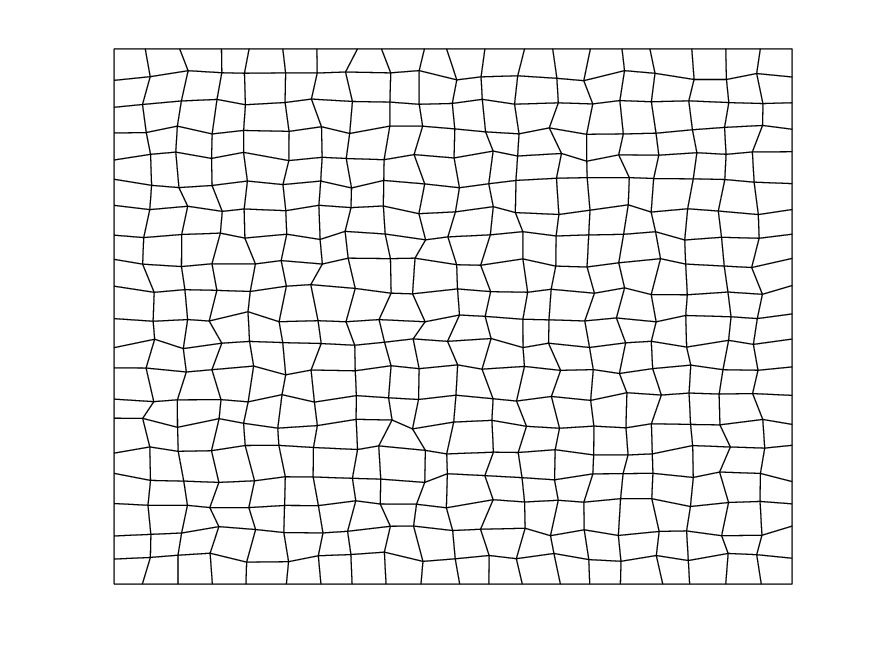}\label{meshrex}}
	    \subfloat[]{\includegraphics[width=0.5\linewidth]{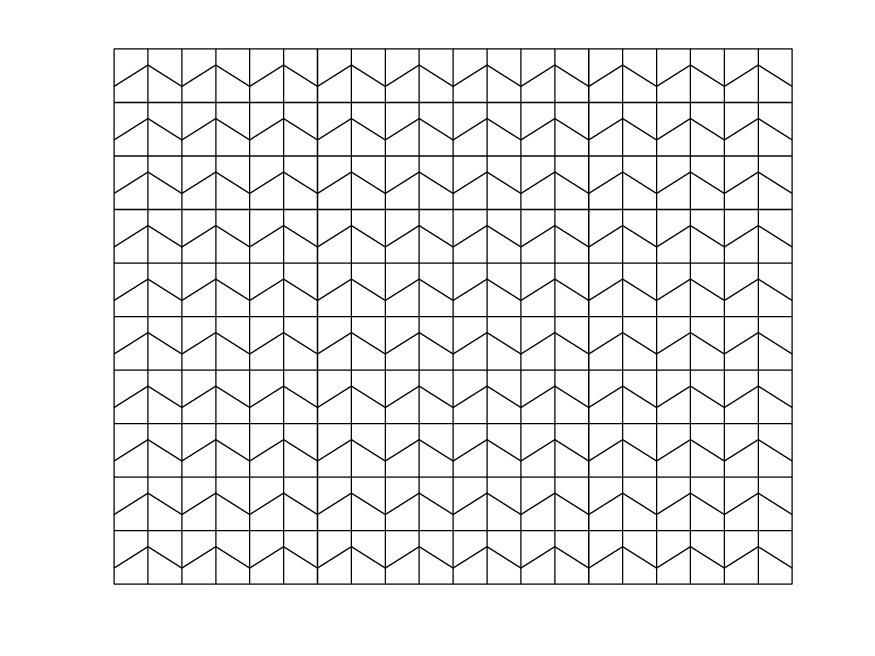}\label{meshtex}}
	    \caption{Two types of $20\times 20$ quadrilateral mesh, both with the degree of distortion $\omega=0.8$. a) Random mesh; b) Trapezoidal mesh.}
	    \label{meshd}
	\end{figure}

    Two types of quadrilateral meshes are considered. The first type, called random mesh, is constructed by adding a non-uniform random perturbation $(\frac{1}{4}r_{i,j}^x h_x,\frac{1}{4}r_{i,j}^y h_y)$ to each node of a $I\times J$ regular mesh with size $h_x\times h_y$. Here $r_{i,j}^x$, $r_{i,j}^y$ are independent random numbers uniformly distributed in $[-\omega,\omega]$, where $\omega\in [0,1]$ is a constant called \textit{the degree of distortion}. At the boundary, $r_{0,j}^x=r_{I,j}^x=0$ for $j=0,\cdots,J$, $r_{i,0}^y=r_{i,J}^y=0$ for $i=0,\cdots,I$. An example of random mesh with $\omega=0.8$ is shown in Figure \ref{meshrex}. 
    Another type is the trapezoidal mesh. It can be constructed by adding a perturbation $(-1)^{i+1}\omega h_y$ on nodes $(x_i,y_{2j-1})$ of the regular mesh, where $\omega\in [0,1]$ is also called \textit{the degree of distortion}. We use $\omega=0.8$ in this example and the meshes are shown in Figure \ref{meshtex}.
    
	Similar as in \cite{HanTwo}, we interpolate the solution on coarse quadrilateral meshes by \eqref{gs2d1} to get the solution on the nodes of a $400\times 400$ regular mesh. Then the extended $l^2$ error is defined by $$\frac{1}{401}\sqrt{\sum_{m=1}^{\tM}\sum_{i=0}^{400}\sum_{j=0}^{400}\left(\Psi_{m}^{interp}(x_i,y_j)-\Psi_{m}^{exact}(x_i,y_j)\right)^2},$$
    where $\Psi_{m}^{exact}(x,y)$ is the exact solution on the $400\times 400$ regular mesh. 
    However, for general quadrilateral meshes, the interpolated solution $\bpsi_{i,j}(x,y)$ \eqref{gs2d1} is not always a good approximation inside the whole cell $C_{i,j}$. The four edge centers can determine a rectangular domain (denoted by $C_{i,j}^{valid}$) whose edges are parallel to the $x$ or $y$ axis and pass through the four edge centers. Recall that the basis functions are in the form of \eqref{eq:basefunction}, the interpolated solution is a good approximation inside $C_{i,j}^{valid}$ but not 
    for all $(x,y)\in C_{i,j}$,  
    $\bpsi_{i,j}(x,y)$ may blow up when $\varepsilon \to 0$. 
    For general quadrilateral meshes, there may exist rectangular subdomains that are not contained in any $C_{i,j}^{valid}$. When a node $(x_r,y_r)$ of the $400\times 400$ regular mesh are inside a subdomain $C_r$ that does not belong to any $C_{i,j}^{valid}$, we compute the inflow boundary conditions at the four edge centers of $C_r$ by interpolations inside $C_{i,j}^{valid}$, then interpolate again by the local solution of the form \eqref{gs2d1} inside $C_r$ (The form of \eqref{gs2d1} is for $C_{i,j}$ as in \eqref{Cij}, we need to replace $x_{i-1}$, $x_i$, $y_j$, $y_{j-1}$ as well).
    
    \begin{figure}[!htbp]
	    \centering
	    \subfloat[]{\includegraphics[width=0.5\linewidth]{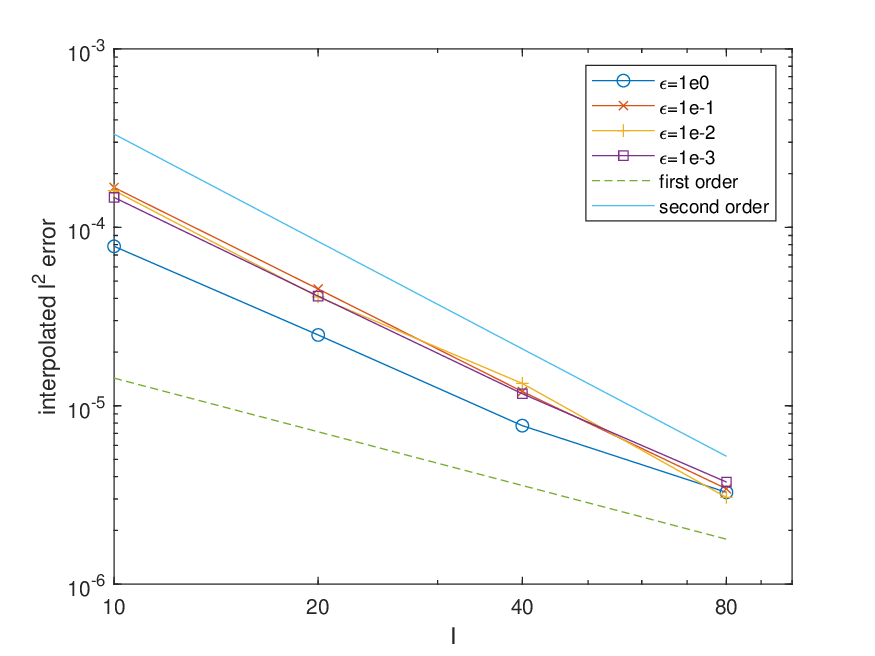}\label{err2drM12}}
	    \subfloat[]{\includegraphics[width=0.5\linewidth]{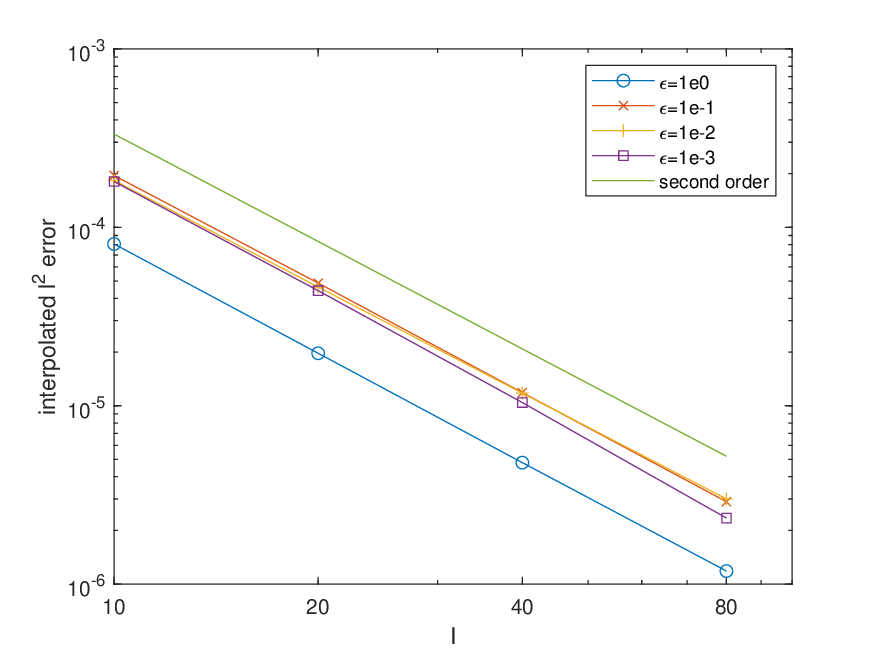}\label{err2dtM12}}
	    \quad
	    \vspace{-0.1cm}
        \subfloat[]{\includegraphics[width=0.5\linewidth]{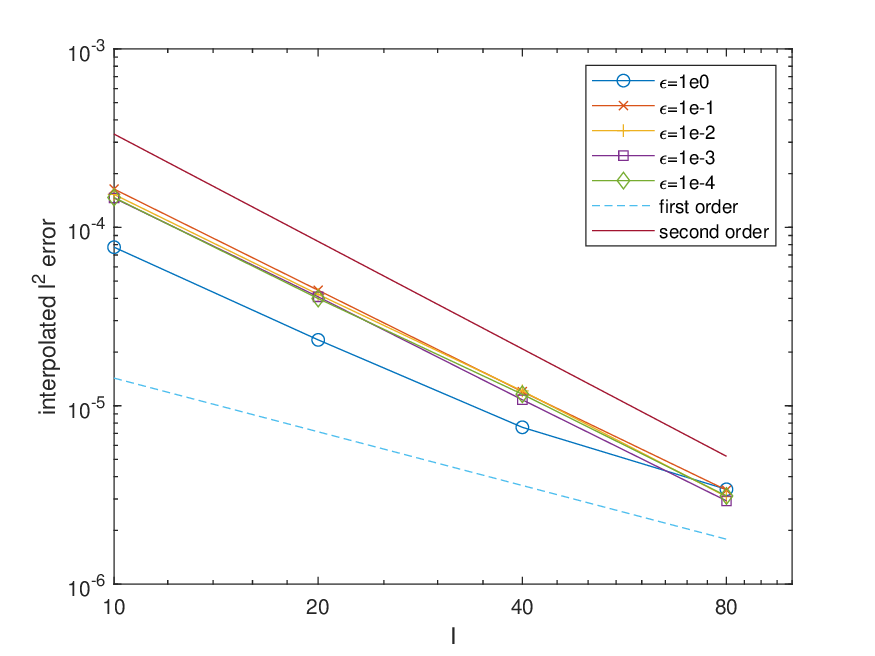}\label{err2dr}}
	    \subfloat[]{\includegraphics[width=0.5\linewidth]{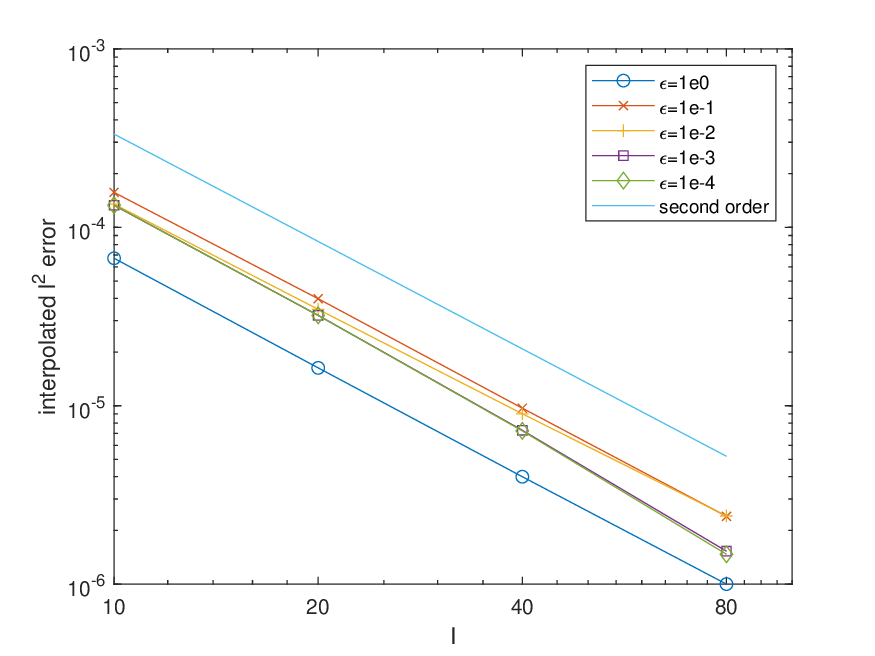}\label{err2dt}}
    \caption{Extended $l^2$ error on general quadrilateral meshes with different $\varepsilon$ and $I$: (a) random meshes with $\tM=12$; (b) trapezoidal meshes $\tM=12$; (c) random meshes with $\tM=24$; (d) trapezoidal meshes $\tM=24$. In all tests, we take the degree of distortion $\omega=0.8$.}
	    \label{err2d}
	\end{figure}
    
    \blue{Extended $l^2$ errors on two types of quadrilateral meshes are plotted in Figure \ref{err2drM12}-\ref{err2dt} for different numbers of spatial cells $\tM$ and rescaled mean free path $\varepsilon$. Uniform second order convergences with respect to $\varepsilon$ are observed for the trapezoidal mesh, while when $\varepsilon=1$ and $I$ becomes bigger, the convergence order reduces to first order using random meshes. But all the convergence results are independent of number of spatial cells $\tM$.}
    
\begin{figure}[!htbp]
		\centering 
		\subfloat[]{\includegraphics[width=.33\textwidth]{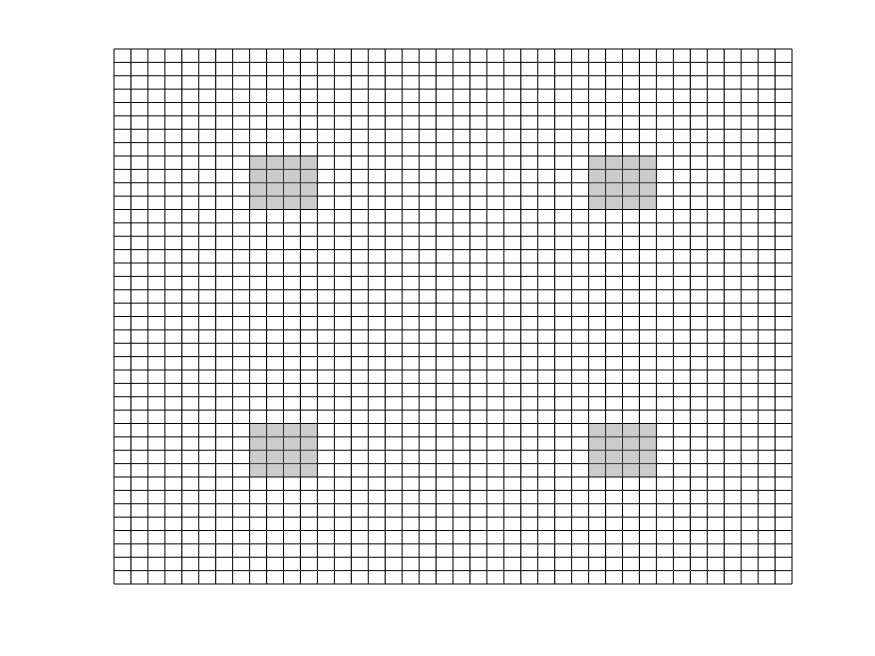}\label{meshu}}
		\subfloat[]{\includegraphics[width=.33\textwidth]{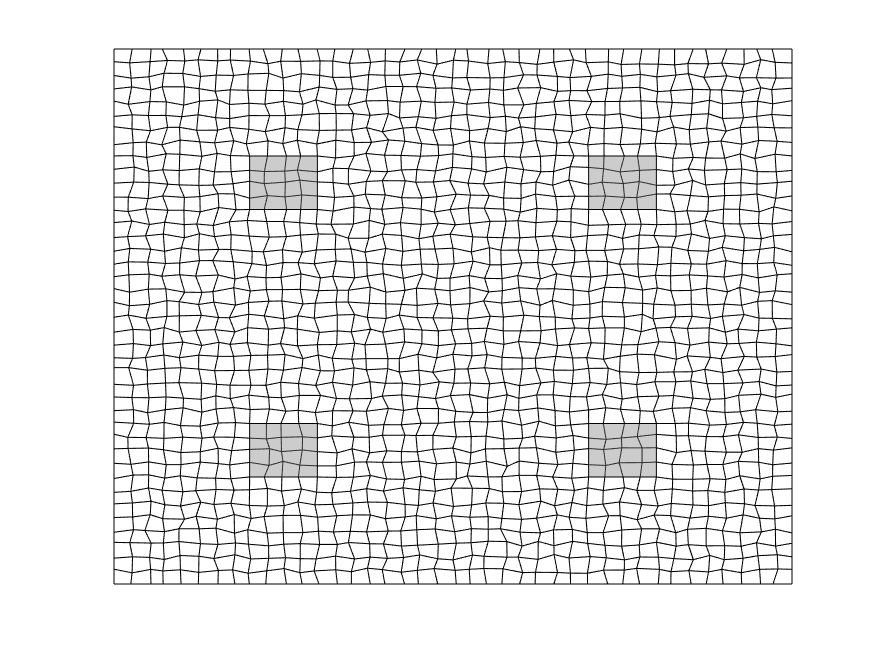}\label{meshr}}
		\subfloat[]{\includegraphics[width=.33\textwidth]{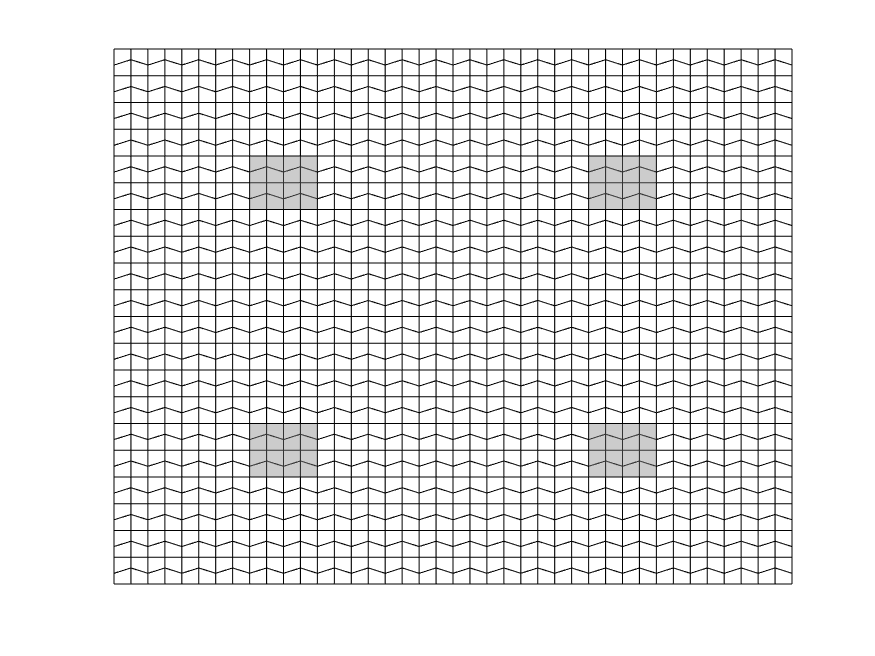}\label{mesht}}\quad 
		\subfloat[]{\includegraphics[width=.33\textwidth]{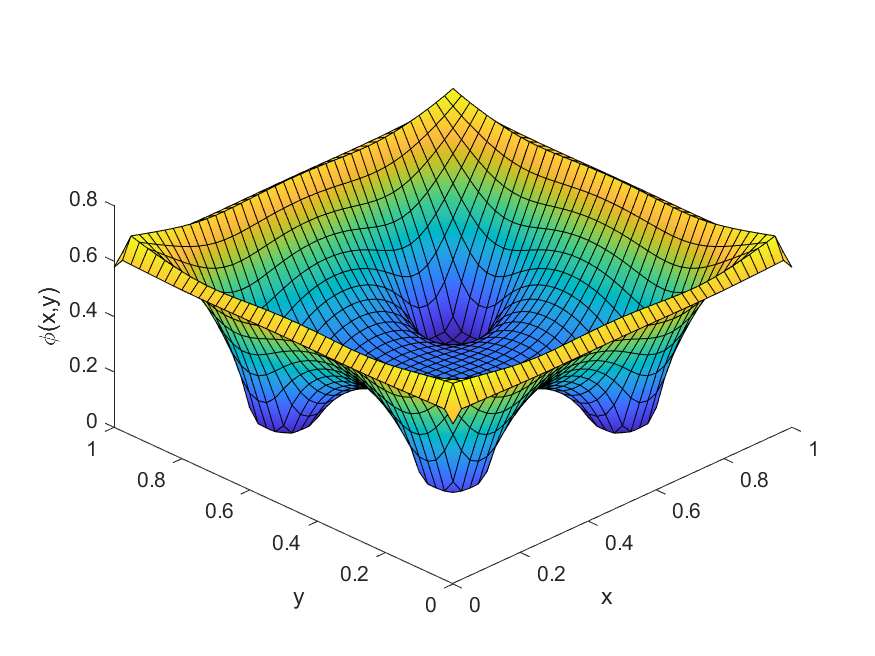}\label{fluxu}}
		\subfloat[]{\includegraphics[width=.33\textwidth]{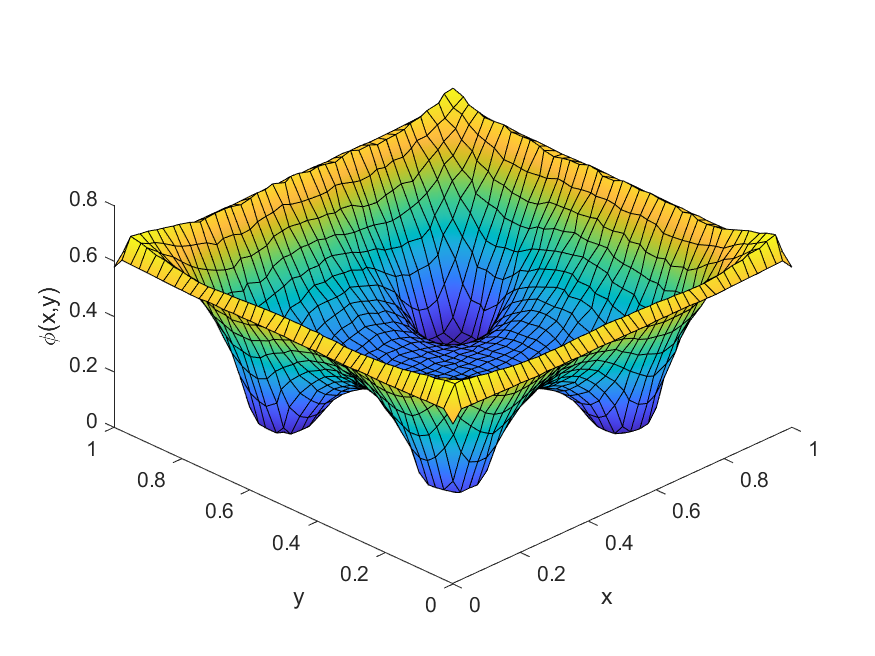}\label{fluxr}}
		\subfloat[]{\includegraphics[width=.33\textwidth]{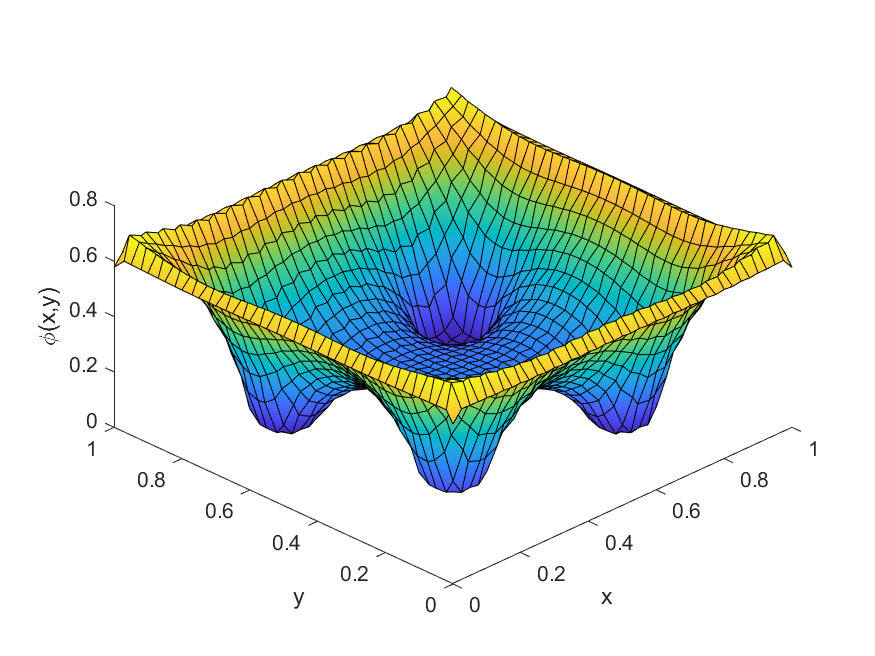}\label{fluxt}}
		\caption{Example 3. The numerical solutions of three different types of meshes with $\tM=24$, $I=J=40$. Meshes are shown in Figure \ref{meshu}-\ref{mesht}, where the gray areas are $\Omega_C$, and in distorted meshes we take the degree of distortion $\omega=0.8$. Density fluxes $\phi(x,y)=\sum_{k\in\bV}\bomega_k\psi_k(x,y)$ are shown in Figure \ref{fluxu}-\ref{fluxr}.}
    \end{figure}

    	\paragraph{Example 3: } In this 2D example, the computational domain consists of both diffusion and transport regions. Four small subdomains are optically thin, and they are surrounded by optical thick materials. The inflow boundary condition is anisotropic. Thus both boundary and interface layers appear. Let
    \begin{equation*}
		\begin{aligned}
			&(x,y)\in\Omega= [0,1]\times[0,1],\quad \Omega_C=\Omega_1\cup \Omega_2 \cup \Omega_3 \cup \Omega_4, \\ &\Omega_1=[0.2,0.3]\times [0.2,0.3],\quad \Omega_2=[0.2,0.3]\times [0.7,0.8],\quad \Omega_3=[0.7,0.8]\times [0.2,0.3],\quad \Omega_4=[0.7,0.8]\times [0.7,0.8]; \\
			&\sigma_T(x,y)=0.2+x+y, \quad \sigma_a(x,y)=0.1+x^2+y^2, \quad \varepsilon=1e-4,\quad (x,y)\in \Omega\setminus \Omega_C; \\
			&\sigma_T(x,y)=1, \quad \sigma_a(x,y)=0.5, \quad \varepsilon=1,\quad (x,y)\in \Omega_C; \\
			&\psi_{l,m}(y_{j-\frac{1}{2}})=c_m,\ c_m>0,\quad \psi_{r,m}(y_{j-\frac{1}{2}})=-c_m,\ c_m<0, \quad j=1,2\cdots J;\\
			&\psi_{b,m}(x_{i-\frac{1}{2}})=s_m,\ s_m>0,\quad \psi_{u,m}(x_{i-\frac{1}{2}})=-s_m,\ s_m<0,  \quad i=1,2\cdots I. \\ 
		\end{aligned}
	\end{equation*}
	The numerical solutions of three different meshes as in Figure \ref{meshu}-\ref{mesht} are displayed in Figure \ref{fluxu}-\ref{fluxr}. Boundary layers can be observed, and the solution behavior highly depends on the material properties inside $\Omega_C$.  

\paragraph{Case I}
We check the dependence of $T_{off}$, $T_{on}$, $T_{gmres}$ on rescaled mean free path $\varepsilon$ (inside $\Omega_C$) and number of spatial cells $I^2$ in Case I, i.e. only the boundary conditions vary. When the number of discrete ordinates is $\tM=24$ and space mesh is $I=J=20$, in Table \ref{2Dtbvarepsilon}, the runtime of offline stage ($T_{off}$), online stage ($T_{on}$), solving large sparse system by preconditioned-GMRES ($T_{gmres}$, for comparison) for different rescaled mean free paths $\varepsilon$ are shown. The computational costs of all stages are independent of $\varepsilon$. In Table \ref{2DItbM}, we fix $\varepsilon=1$, $I=J=20$, and show the runtime of different stages for different $\tM$. Increasing rate of $T_{off}$ w.r.t $\tM$ is about $M^{2.8}$, and increasing rate of $T_{on}$ is about $M^{2}$ and $GMRES$ is $M^{2.5}$. In Table \ref{2DItb}, we fix $\varepsilon=1$, $\tM=24$, and show the runtime of different stages for different $I$ ($J=I$). Increasing rates of $T_{off}$, $T_{on}$ and $T_{gmres}$ w.r.t $I$ are respectively around $I^{3.8}$, $I^{3}$ and $I^{3.3}$.
\begin{table}[!htbp]
    \centering
    \subfloat[]{
            \begin{tabular}{|c|c|c|c|}
				\hline
				$\varepsilon$ & $T_{off}$ & $T_{on}$ & $T_{gmres}$  \\
				\hline
				1 & 18.44 & 0.19 & 2.40\\
				\hline
				$10^{-1}$ & 19.51& 0.23 & 2.75\\
				\hline
				$10^{-2}$ & 19.39& 0.19 & 2.85\\
				\hline
				$10^{-3}$ & 18.42& 0.20 & 2.70\\
				\hline
				\end{tabular}
            \label{2Dtbvarepsilon}}
            \hspace{0.5cm}
        \subfloat[]{		\begin{tabular}{|c|c|c|c|}
		    \hline
			$\tM$ & $T_{off}$ & $T_{on}$ & $T_{gmres}$\\
			\hline
			24 & 18.44 & 0.19& 2.40\\
			\hline
			40 & 83.29& 0.77& 8.52\\
			\hline
			60 & 245.16& 1.74& 22.25\\
			\hline
			84 & 591.17 & 2.42 & 53.46\\
			\hline
		\end{tabular}\label{2DItbM}}
		\hspace{0.5cm}
		\subfloat[]{\begin{tabular}{|c|c|c|c|}
		    \hline
		    $I^2$ & $T_{off}$ & $T_{on}$ & $T_{gmres}$ \\
			\hline
			100 & 1.15 & 0.03& 0.26 \\
			\hline
			400 & 18.44 & 0.19& 2.40 \\
			\hline
			1600 & 336.11& 1.55& 24.27 \\
			\hline
			6400 & 2925.88 & 15.53 & 239.09 \\
			\hline
		\end{tabular}
		\label{2DItbI}}
    \caption{Example 3. 2D Case I. Run time (seconds) of different stages in the fast solver with different rescaled mean free path $\varepsilon$, number of discrete ordinates $\tM$ and number of spatial cells $I^2$.}
    \label{2DItb}
\end{table}

\paragraph{Case II}
 We check the dependence of $T_{off}$, $T_{on}$, $T_{gmres}$ on $\tM$ and number of spatial cells $I^2$ in Case II. $\sigma_T$, $\sigma_a$ and $\varepsilon$ may vary inside $\Omega_C$ which is a small domain compared to $\Omega$. The number of spatial cells inside $\Omega_C$ is taken to be $4I^2/100$, so that the mesh sizes inside and outside of $\Omega_C$ are the same. We fix $I=J=20$ and record $T_{off}$, $T_{on}$, $T_{gmres}$ for different numbers of discrete ordinates $\tM$=24, 40, 60, 84. As shown in Table \ref{2DIItbM}, increasing rate of $T_{off}$ w.r.t $\tM$ is about $M^{2.7}$, and increasing rate of $T_{on}$ is about $M^{1.4}$ and $T_{gmres}$ is $M^{2.5}$. On the other hand, we fix $\tM=24$ and record $T_{off}$, $T_{on}$, $T_{gmres}$ for different numbers of spatial cells $I^2$=100, 400, 1600, 6400. As shown in Table  \ref{2DIItbI}, increasing rates of $T_{off}$, $T_{on}$ and $T_{gmres}$ w.r.t $I$ are respectively around $I^{3.6}$, $I^{2.5}$ and $I^{3.3}$. Run time of the online stage is much faster. When $\Omega_C$ is smaller, our online stage can be even cheaper. Moreover, our on line stage can be easily accelerated by parallelization for different spacial cells.

\begin{table}[!htbp]
	\centering
	\subfloat[]{\begin{tabular}{|c|c|c|c|}
		\hline
		$\tM$ & $T_{off}$ & $T_{on}$ & $T_{gmres}$ \\
		\hline
		24 & 34.72& 0.78 & 2.37 \\
		\hline
		40 & 128.98& 1.30 & 8.16\\
		\hline
		60 & 447.69& 2.71 & 22.49 \\
		\hline
		84 & 1035.95 & 4.17 & 51.38 \\
		\hline
	\end{tabular}\label{2DIItbM}}
	\hspace{0.5cm}
	\subfloat[]{\begin{tabular}{|c|c|c|c|}
		\hline
		$I^2$ & $T_{off}$ & $T_{on}$ & $T_{gmres}$ \\
		\hline
		100 & 2.67& 0.21 & 0.26 \\
		\hline
		400 & 34.72& 0.78& 2.37\\
		\hline
		1600 & 627.41& 4.61 & 24.80 \\
		\hline
		6400 & 4352.87 & 37.19 & 237.95 \\
		\hline
	\end{tabular}\label{2DIItbI}}
\caption{Example 3. 2D Case II. Run time (seconds) of different stages in the fast solver (offline/online) with different numbers of discrete ordinates $\tM$ and numbers of spatial cells $I^2$.}
\label{2DIItb}
\end{table}


	\section{Conclusion and discussion}

	
	This paper presents a general approach to decomposing TFPM of RTE into offline/online stages in the two cases illustrated in the introduction. In Case I, only the right-hand side of the linear system $A\balpha=\bb$ varies. In Case II, not only the right-hand side $\bb$ but also the coefficient matrix $A$ changes. The expensive offline stage is only calculated once, and the cheap online stage is updated for each different parameter chosen from a large data set.

	Our scheme can be understood as a preconditioner that transforms the coefficient matrix into a block-diagonal form. Most expensive operations like matrix-matrix multiplication or $QR/LU$ decomposition are finished at the offline stage.  One may try to inverse the coefficient matrix at the offline stage directly and then use matrix-vector multiplication to get the solution at the online stage. However, it is hard to find an efficient way to inverse the coefficient matrix directly, especially in 2D. Other fast solvers can be understood as a preconditioner that can deal with multiple right-hand sides, including Block GMRES\cite{simoncini1996convergence}, Block BiCG-STAB \cite{el2003block}, etc. For example, in 1D, the linear system is block tri-diagonal, and classical Block LU decomposition can be used. At the offline stage, the coefficient matrix $A$ is decomposed into a lower bi-diagonal block matrix $L$ and an upper bi-diagonal block matrix $R$. At the online stage, two bi-diagonal systems are solved for different $\bb$. In our solver, when $\bb$ changes, some matrix-vector multiplications are needed to update the right-hand side of small local systems $\bb_i$, and then some small upper/lower triangular systems have to be solved. The cost at the online stage in Block LU is comparable to our solver, but it is not cheap in 2D since the sparsity of $A$ is different. In \cite{chen2019low}, a fast low-rank method has been developed for linear RTE, which can be considered as an offline/online algorithm as well. The low-rank structure is regarded at the offline stage, and the cost at the online stage is much lower. There exist other fast solvers or more advanced iteration methods for linear RTE \cite{ren2019fast,fan2019fast}. However, their performances for problems that exhibit boundary or interface layers are barely considered. 
	Compared with the above-mentioned methods, our approach has three benefits: 1) Since the coefficient matrix has been transformed into a block-diagonal form, solutions at different regions can be solved in parallel at the online stage; 2) 
when $\sigma_a$, $\sigma_S$ change locally, the cost to update the preconditioner depends only on the region size where $\sigma_a$, $\sigma_S$ vary; 3) Our approach is applicable to problems when the material optical properties vary a lot in the computational domain.

In the numerical tests, the restarted GMRES solver with block-diagonal right-preconditioner and ILU right-preconditioner is a standard solver in numerical algebra. 
We observe that the iteration steps of restarted GMRES solver may increase dramatically for the 2D problem with random meshes when the solution exhibits boundary and interface layers. Sometimes it does not converge. This phenomenon can not be observed with only optical thin or thick media. This indicates that the performances of iterative solvers depend on the properties of the materials and the meshes.

\blue{The operation is applicable in 3D, the bottle-neck is the storage. In 3D, more velocity coordinates are needed and the operation matrix is much larger. The current simulations are all done in an Inter Xeon Processor (Skylake, IBRS) @ 2.39 GHz, which is not enough for 3D simulations. Besides, it is important to find an efficient way to do $QR$ decomposition and store all necessary information when the dimension becomes higher.} To solve inverse RTE, one particularly interesting case is how to efficiently update the solutions when $\sigma_a$, $\sigma_T$ change a little bit. It should be noted that when cross-sections vary in whole area $\Omega$, all the non-zero elements in $A$ change, and our solver is no longer fast. We will investigate this case in our future work.
	

	\begin{appendices}
		\renewcommand{\theequation}{A.\arabic{equation}} 
		\section{Quadrature choices in discrete-ordinate methods}
		\subsection{1D case}
		We choose a quadrature set of size $2M$: $\{\mu_m,\omega_m|m \in V \}$, where $V$ is the order set $\{-M,-M+1,\cdots,-2,-1,1,2,\cdots,\\M-1,M\}$,
		and weights $\omega_m$ are normalized by
		\eq$\sum_{k \in V}\omega_{k}=1.\label{cons1}$
		
		By classical asymptotic analysis in \cite{larsen1987asymptotic}, if the quadrature set satisfies:
		\eq$\sum_{k \in V}\omega_{k}\mu_{k}=0,\quad \sum_{k \in V}\omega_{k}\mu_{k}^2=\frac{1}{3}, \label{cons2}$
		then $\psi_{m}=\phi+\mathcal{O}(\varepsilon)$ when $\varepsilon\to 0$, where $\phi$ is the solution for diffusion limit equation \eqref{1ddiff}.
		
		Specifically, we take Gauss-Legendre quadrature: $\{\mu_m|m\in V\}$ is $2M$ distinct roots of Legendre polynomials $P_{2M}(x)$ with degree $2M$, ordered as
		\eq$-1<\mu_{-M}<\mu_{-M+1}<\cdots<\mu_{-2}<\mu_{-1}<0<\mu_1<\mu_2<\cdots<\mu_{M-1}<\mu_{M}<1,$
		and 
		\eq$\omega_m=\frac{2}{(1-\mu_m^2)[P_{2M}'(\mu_m)]^2}.$
		
		Gauss-Legendre quadrature satisfies constraints \eqref{cons1} and \eqref{cons2} mentioned above, and guarantees the symmetry of $\omega_{m}$ and $\mu_m$:
		\eq$\omega_m=\omega_{-m},\quad \mu_{m}=\mu_{-m}.$
		\subsection{2D case}
		The Gaussian quadrature set $S_N$ for $\theta \in [0,\frac{\pi}{2}]$ is generated in the following way:
		\begin{itemize}
			\item each quadrant has $M=N(N+1)/2$ ordinates.
			\item each quadrant has $N$ distinct $\zeta_n$, which are the positive roots of Legendre polynomials $P_{2N}$, ordered as $0<\zeta_1<\zeta_2<\cdots<\zeta_N<1$.
			\item each $\zeta_n$ correspond to $m$ distinct $\theta_{n,i}=\frac{2i-1}{4n}\pi$, $i=1,2,\cdots,m$ and the same weights   
			\eq$\bomega_n=\frac{1}{n(1-\zeta_n^2)[P_{2N}'(\mu_i)]^2}.$
			\item reorder $\{(\theta_m,\bomega_m,\zeta_m)|m=1,2,\cdots,M\}$ by 
			$$\{(\theta_{1,1},\bomega_1,\zeta_1),(\theta_{2,1},\bomega_2,\zeta_2),(\theta_{2,2},\bomega_2,\zeta_2),(\theta_{3,1},\bomega_3,\zeta_3),\cdots,(\theta_{N,N},\bomega_N,\zeta_N)\}$$
		\end{itemize}
		
		Then the remainder of the quadrature set can be constructed by symmetry:
		\eq$
		\begin{aligned}
		&\theta_m=\theta_{m+M}-\frac{\pi}{2}=\theta_{m+2M}-\pi=\theta_{m+4M}-\frac{3}{2}\pi, \\
		&\bomega_m=\bomega_{m+M}=\bomega_{m+2M}=\bomega_{m+4M}, \\
		&\zeta_m=\zeta_{m+M}=\zeta_{m+2M}=\zeta_{m+4M}, \\
		\end{aligned}$
		for $m=1,2,\cdots,M$ and 
		\eq$c_m=(1-\zeta_m^2)^{\frac{1}{2}}\cos\theta_m,\quad s_m=(1-\zeta_m^2)^{\frac{1}{2}}\sin\theta_m, \quad m\in \bV.$
		
		The generated Gaussian quadrature set $S_N$ guarantees  
		\eq$\sum_{k \in \bV}\bomega_k=1,\quad \sum_{k \in \bV}\bomega_k c_k=0,\quad \sum_{k \in \bV}\bomega_k s_k=0,\quad \sum_{k \in \bV}\bomega_k c_k s_k=0,\quad \sum_{k \in \bV}\bomega_k (c_k^2+s_k^2)=\frac{2}{3},$
		which indicates the discrete-ordinate equations \eqref{eqad2} and \eqref{2deq} converges to the same diffusion limit when $\varepsilon\to 0$ \cite{tang2009uniform}.
		\section{Explicit form of block tri-diagonal system in 1D}
		\renewcommand{\theequation}{B.\arabic{equation}} 
		\begin{equation}
		A=\left(\begin{array}{c;{2pt/2pt}c;{2pt/2pt}c;{2pt/2pt}c;{2pt/2pt}c}
		A_1^r(x_0) & & & &\\
		A_1^l(x_1) & -A_2^l(x_1) & & &\\\hdashline[2pt/2pt]
		-A_1^r(x_1) & A_2^r(x_1) & & &\\
		& A_2^l(x_2) & -A_3^l(x_2) & &\\\hdashline[2pt/2pt]
		& -A_2^r(x_2) & A_3^r(x_2) & &\\ 
		& & A_3^l(x_3) & -A_4^l(x_3) &\\\hdashline[2pt/2pt]
		& & \ddots& \ddots & \ddots \\\hdashline[2pt/2pt]
		& & & -A_{I-1}^r(x_{I-1})& A_{I}^r(x_{I-1}) \\
		& & & & A_I^l(x_I) \\
		\end{array}\right)
		\label{defA}
		\end{equation}
    \eq$b=\left(\begin{array}{c}
		\bpsi_l^b \\
		\bm{0}_{M/2}\\\hdashline[2pt/2pt]
		\bm{0}_{M/2}\\
		\bm{0}_{M/2}\\\hdashline[2pt/2pt]
		\bm{0}_{M/2}\\
		\bm{0}_{M/2}\\\hdashline[2pt/2pt]
		\vdots \\\hdashline[2pt/2pt]
		\bm{0}_{M/2}\\
		\bpsi_r^t\\
		\end{array}\right)$
 	\end{appendices}
 
 \section*{Acknowledgments}
 Min Tang is partially supported by Shanghai Pilot Innovation project, 21JC1403500, the Strategic Priority Research Program of Chinese Academy of Sciences, XDA25010401 and NSFC12031013. Jingyi Fu is partially supported by Shanghai Pilot Innovation project, 21JC1403500 and would like to thank Yihong Wang for providing the generation of distorted meshes.
 
	\renewcommand\refname{Reference}
    \bibliography{AM}
 	\bibliographystyle{CICP}


\end{document}